\documentclass[11pt]{article}

\usepackage{amssymb,amsmath,euscript,bbm,xcolor,color,graphicx,epstopdf,bm}
\usepackage{verbatim} 
\newtheorem{proposition}{Proposition}[section]
\newtheorem{lemma}[proposition]{Lemma}
\newtheorem{theorem}[proposition]{Theorem}

\def\ep{\varepsilon}

\def\l{{\langle}}
\def\r{\rangle}

\newcommand{\wt}{\widetilde}
\def\R{{\mathbb R}}
\def\S{{\mathbb S}}

\def\E{{\mathbb E}}
\def\P{{\mathbb P}}

\newcommand{\prooftheo}[1]{ \textsc{Proof of Theorem} \ref{#1} }

\newcommand{\pk}[1]{\mathbb{P}\left\{{#1}\right\}}
\newcommand{\QED}{\hfill $\Box$}

\newcommand{\COM}[1]{}

\def\IF{\infty}

\def\bqn#1{ { \begin{eqnarray} #1 \end{eqnarray}}}

\newcommand{\BS}{\begin{sat}}
	\newcommand{\ES}{\end{sat}}
\newcommand{\BT}{\begin{theorem}}
	\newcommand{\ET}{\end{theorem}}
\newcommand{\BK}{\begin{korr}}
	\newcommand{\EK}{\end{korr}}

\newcommand{\BEX}{\begin{example}}
	\newcommand{\EEX}{\end{example}}

\newcommand{\BD}{\begin{de}}
	\newcommand{\ED}{\end{de}}

\newcommand{\BIT}{\begin{itemize}}
	\newcommand{\EIT}{\end{itemize}}
\newcommand{\BDI}{\begin{description}}
	\newcommand{\EDI}{\end{description}}

\newcommand{\BRM}{\begin{remark}}
	\newcommand{\ERM}{\end{remark}}

\newcommand{\BEL}{\begin{lemma}}
	\newcommand{\EEL}{\end{lemma}}

\newcommand{\BQN}{\begin{eqnarray}}
\newcommand{\EQN}{\end{eqnarray}}
\newcommand{\BQNY}{\begin{eqnarray*}}
	\newcommand{\EQNY}{\end{eqnarray*}}

\makeatletter \@addtoreset{equation}{section} \makeatother
\newcommand {\qed}%
{%
    {}\hfill
    {}\hfill
    {$\square $}%
    \vspace {0.3cm}%
    \pagebreak [2]%
    \par
}%
\newenvironment{proof}[1]{%
    \vspace{0.3cm}%
    \pagebreak [2]%
    \par%
    \noindent {\bf  Proof.~#1\ }}{\qed}%
\newenvironment{example}{%
    \vspace{0.3cm} \pagebreak [2]%
    \par%
    \refstepcounter{proposition}%
    \noindent%
    {\bf  Example~\theproposition\ }}{}%
\newenvironment{remark}{%
    \vspace{0.3cm} \pagebreak [2]%
    \par%
    \refstepcounter{proposition}
    \noindent%
    {\bf Remark~\theproposition\  }}{}%
\begin{document}

\title {Extremes of Spherical Fractional Brownian Motion}
\author{Dan Cheng \thanks{Supported by NSF Grant DMS-1811632.} \\ Arizona State University
 \and Peng Liu \thanks{Supported by Swiss National Science Foundation Grant 200021-175752/1.} \\ University of Lausanne}

\maketitle

\begin{abstract}
Let $\{B_\beta (x), x \in \mathbb \S^N\}$ be a fractional Brownian motion on the $N$-dimensional unit sphere $\S^N$ with Hurst index $\beta$. We study the excursion probability $\pk{\sup_{x\in T} B_\beta(x) > u }$ and obtain the asymptotics as $u\to \infty$, where $T$ can be the entire sphere $\S^N$ or a geodesic disc on $\S^N$.
\end{abstract}



\section{Introduction}
Let $\{X(t), t\in T \}$ be a real-valued Gaussian random field living on some parameter space $T$. The extremes, especially excursion probabilities $\P\{\sup_{t\in T} X(t)> u \}$, of the field have been extensively studied in the literature due to the importance in both probability theory \cite{Rice:1945,ChanL06,Debicki:2016,HP99,Liu12,Sun93} and statistical applications such as the p-value computation for controlling the family-wise error \cite{TaylorW07,TaylorW08}, nonparametric density estimation \cite{Bickel:1973,Huckemann:2016,Qiao:2017} and construction of confidence bands \cite{MaYang:2012,WangYang:2007}. We refer to the survey \cite{Adler00} and monographs \cite{Piterbarg:1996,Adler07,AzaisBook09} for the history, recent developments and more related applications on this subject.

Recently, the study of random fields on spheres is attracting more and more attention due to vast applications in astronomy \cite{MP11}, spatial statistics \cite{Gelfand:2010,Stein:1999}, geoscience \cite{OhLi04,Jeong:2017} and environmental sciences \cite{Stein:2007}. In particular, Istas \cite{Istas:2005,Istas:2006} introduced spherical fractional Brownian motion (abbreviated as SFBM throughout this paper) on spheres and studied the Karhunen-Lo\`eve expansion and other properties. As an important extension to the classical fractional Brownian motion on Euclidean space, it would be very useful and valuable to study the excursion probability of SFBM, which is the main purpose of this paper.

Let $o$ be a fixed point on the $N$-dimensional unit sphere $\mathbb{S}^N \subset \R^{N+1}$. The SFBM, denoted by $B_\beta= \{B_\beta (x), x \in \mathbb S^N\}$, is defined in Istas \cite{Istas:2005} as a centered real-valued Gaussian random field on $\S^N$ such that $B_\beta(o)=0$ and
\begin{equation}\label{eq:SFBM-var}
\E[B_\beta(x)-B_\beta(y)]^2= d^{2\beta}(x,y), \quad \forall x,y\in \mathbb{S}^N,
\end{equation}
where the index $\beta\in (0, 1/2]$ and  $d(\cdot,\cdot)$ is the spherical distance on $\mathbb{S}^N$, that is $d(x,y)= \arccos \l x, y \r$, $\forall x, y\in \mathbb{S}^N$. Here $\l \cdot, \cdot \r$ is the usual inner product in $\R^{N+1}$. It follows immediately that the covariance structure is given by
\begin{equation}\label{eq:SFBM-cov}
{\rm Cov}(B_\beta(x), B_\beta(y))=\frac{1}{2}\big(d^{2\beta}(x,o)+d^{2\beta}(y,o)-d^{2\beta}(x,y) \big).
\end{equation}

In this paper, we shall study the asymptotics of the excursion probability $\P\{\sup_{x\in T} B_\beta(x)> u \}$ as $u\to \infty$. Two cases for the parameter set $T$ are considered separately: (i) $T=\S^N$ and (ii) $T=T_a:=\{x\in\mathbb{S}^N: d(x,o)\le a\}$, where $a\in (0, \pi)$. In other words, $T_a$ is the geodesic disc on $\S^N$ of radius $a$ centered at $o$, so that $T_a=\S^N$ when $a=\pi$. Notice that, the maximum of the variance function of $B_\beta(x)$ over $T$ will be attained at a single point for case (i) and on the boundary set $\{x\in\mathbb{S}^N: d(x,o)= a\}$, which is in fact an $(N-1)$-dimensional sphere, for case (ii), respectively, making the latter case more challenging.

Since the sphere $\S^N$ is not an Euclidean space, it would be hard to apply directly the traditional double sum method over $\S^N$ to derive the asymptotics of the excursion probability. Instead, we shall apply the main technique in Cheng and Xiao \cite{ChengXiaoSphere} to consider the SFBM as a Gaussian random field on Euclidean space by using spherical coordinate transformation. In such way, we can study the local behaviors of the standard deviation and correlation functions of the field under spherical coordinates (see Lemmas \ref{lemA} and \ref{lemB} below), and then apply the results in Euclidean space (see Lemma \ref{Piterbarg} and Theorem \ref{extension} below) to derive the desired asymptotics of the excursion probabilities in Theorems \ref{Thm-Sphere} and \ref{ThA}. In particular, for case (ii), the maximum of variance is attained on a set of dimension at least one (when $N\ge 2$) and there is no known asymptotic result in the literature except for the two-dimensional case studied in \cite{Debicki:2016}. In order to obtain the asymptotics on a geodesic disc in Theorem \ref{ThA}, we establish an asymptotic result on Euclidean space in Theorem \ref{extension} which is valuable itself in extreme value theory and will have further applications in the future.

The paper is organized as follows. We first introduce the preliminaries, such as spherical coordinate transformation and the Pickands and Piterbarg constants, in Section \ref{sec:prel}; and then study the asymptotics of the excursion probabilities of $B_\beta(x)$ on the entire sphere $\S^N$ and on a geodesic disc in Sections \ref{sec:sphere} and \ref{sec:disc}, respectively. Finally, the Appendix contains some auxiliary results and the proof of Theorem \ref{extension}.


\COM{ Define next
	\begin{equation*}
	\begin{split}
	\sigma(x)&=\sqrt{{\rm Var}(B_\beta(x))}, \quad \tilde{\sigma}(\theta)=\sigma(x)\\
	r(x,y)&={\rm Cov}(B_\beta(x)/d^\beta(x,o), B_\beta(y)/d^\beta(y,o)), \quad \tilde{r}(\theta,\varphi)=r(x,y).\\
	\end{split}
	\end{equation*}
}

\section{Preliminaries}\label{sec:prel}
\subsection{Spherical Coordinates and Notations}\label{subsection:SphereCoor}
For $x=(x_1, \ldots, x_N, x_{N+1})\in \mathbb{S}^N$, its corresponding spherical coordinate $\theta=(\theta_1, \ldots, \theta_N)$ is defined by the following way.
\begin{equation}\label{Eq:spherical coordinate}
x_1 = \cos \theta_1, \quad \ldots, \quad x_{N} = \left(\prod_{i=1}^{N-1}\sin \theta_i \right) \cos \theta_{N}, \quad x_{N+1} = \prod_{i=1}^{N}\sin \theta_i,
\end{equation}
where $\theta\in \Theta:=[0, \pi]^{N-1}\times[0, 2 \pi)$.

Throughout this paper, for two points $x=(x_1, \ldots, x_{N+1})$ and $y=(y_1, \ldots, y_{N+1})$ on $\mathbb{S}^N$, we always denote by $\theta=(\theta_1, \ldots, \theta_N)$ and $\varphi=(\varphi_1, \ldots, \varphi_N)$ the spherical coordinates of $x$ and $y$, respectively. For functions $f(x)$ and $g(x,y)$, $x, y\in \S^N$, we denote by $\tilde{f}(\theta):=f(x)$ and $\tilde{g}(\theta,\varphi):=g(x,y)$ the corresponding functions of $f(x)$ and $g(x,y)$ under spherical coordinates, respectively.

Let $d(\cdot,\cdot)$ denote the spherical distance on $\mathbb{S}^N$ and let $\|\cdot\|$ be the Euclidean norm in $\R^{N+1}$ or in $\R^N$, which will be clear from the context. For a set $D\subset \R^N$, denote by ${\rm mes} (D)$ the measure (volume) of $D$. Denote by $\Psi(u)$ the tail probability of standard normal distribution, that is $\Psi(u)=(2\pi)^{-1/2}\int_u^\infty e^{-v^2/2}dv$. For any two real-valued functions $h_1(u)$ and $h_2(u)$, we say $h_1(u) \sim h_2(u)$ as $u\to u_0\in [-\infty, +\infty]$ if $\lim_{u\to u_0} h_1(u)/h_2(u) =1$.


\subsection{Pickands and Piterbarg Constants}
Let $\{\chi_H(t)$, $t\in \mathbb{R}^N\}$, $H\in (0,1]$, be a  Gaussian random field with mean function
$$\mathbb{E}\left(\chi_H(t)\right)=-\|t\|^{2H}, \quad t\in \R^N,$$
and covariance
$${\rm Cov}(\chi_H(t), \chi_H(s))=\|t\|^{2H}+\|s\|^{2H}-\|t-s\|^{2H}, \quad t,s\in \R^N.$$
Let $\mathcal{H}_{2H}(E)=\mathbb{E}\left\{ \sup_{t\in E }\exp[\chi_H(t)]\right\}$, where $E\subset \R^N$ is a compact set. The \textit{Pickands constant} \cite{Piterbarg:1996} is defined by
\begin{equation}\label{eq:Pickands-Const}
\mathcal{H}_{2H}^N:=\lim_{S\to\infty}\frac{\mathcal{H}_{2H}^N([0,S])}{S},\quad \text{where }  \mathcal{H}_{2H}^N([0,S])=\lim_{S_1\to\infty}\frac{\mathcal{H}_{2H}([0,S]\times[0,S_1]^{N-1})}{S_1^{N-1}}.
\end{equation}
The \textit{Piterbarg constant} \cite{Bai:2018, Piterbarg:1996} is defined by
\begin{equation}\label{eq:Piterbarg-Const}
\mathcal{P}_{2H}^g :=\lim_{S\to\infty}\mathcal{P}_{2H}^g([-S,S]^N), \quad \text{where }\mathcal{P}_{2H}^g(E)=\mathbb{E}\left\{ \sup_{t\in E }e^{\chi_H(t)-g(t)}\right\},
\end{equation}
$E\subset \R^N$ is a compact set and $g$ is a continuous function over $\mathbb{R}^N$. Moreover, let
\begin{equation}\label{eq:M-Const}
\begin{split}
\mathcal{M}_{2H}^g &:=\lim_{S\to\infty}\lim_{S_1\to\infty}\frac{\mathcal{P}_{2H}^g([-S,S]\times[0,S_1]^{N-1})}{S_1^{N-1}},\\
\widehat{\mathcal{M}}_{2H}^g &:=\lim_{S\to\infty}\mathcal{M}_{2H}^g([0,S]), \quad \text{where } \mathcal{M}_{2H}^g([0,S])=\lim_{S_1\to\infty} \frac{\mathcal{P}_{2H}^g([0,S]\times[0,S_1]^{N-1})}{S_1^{N-1}},
\end{split}
\end{equation}
if the limits above exist.

\COM{\subsection{Excursion Probability over the Whole Sphere}
	Once we find the local structure of $\tilde{\sigma}(\theta)$ and $\tilde{r}(\theta,\varphi)$ near $\theta_0$ [see (\ref{Eq:sd sphere}) and (\ref{Eq:corr sphere})], then we can apply Piterbarg's results to get the approximation for $\P\{\sup_{x\in \mathbb{S}^N} B_\beta(x) \ge u \}=\P\{\sup_{\theta\in [0, \pi]^{N-1}\times[0, 2 \pi)} \wt{B}_\beta(\theta) \ge u \}$.
}

\section{Excursion Probability on $\S^N$}\label{sec:sphere}
We first study the excursion probability of $B_\beta$ over the entire sphere. Recall the notations introduced in Section \ref{subsection:SphereCoor}, we have that $\pk{ \sup_{x\in \mathbb{S}^N} B_\beta(x)>u}$ is equivalent to $\P\{ \sup_{\theta\in \Theta} \widetilde{B}_\beta(\theta)>u\}$, where $\wt{B}_\beta(\theta) := B_\beta(x)$ and $\Theta=[0, \pi]^{N-1}\times[0, 2 \pi)$ is an $N$-dimensional rectangle on $\R^N$. Therefore, to establish the asymptotics for the excursion probability of $B_\beta$, we will study the properties of the standard deviation and correlation functions of $\widetilde{B}_\beta$, which is a Gaussian random field living on $\Theta$, and then apply results on extremes for Gaussian random fields on Euclidean space.

In this section, we assume without loss of generality that $B_\beta$ starts at $o=(0,0,\ldots,1,0)\in \R^{N+1}$ whose spherical coordinate is given by $(\pi/2,\ldots,\pi/2,0)\in \Theta \subset \R^N$ according to \eqref{Eq:spherical coordinate}. Denote by $\sigma(x)$ the standard deviation function of $B_\beta(x)$. By \eqref{eq:SFBM-var},
\bqn{ \label{sigmax0}
	\sigma(x)=d^\beta(x,o)= \arccos^{\beta} \l x, o \r = \arccos^{\beta} (x_N), \quad x\in \S^N,}
which attains its unique maximum $\pi^\beta$ at $p:=(0,0,\ldots,-1,0)\in \R^{N+1}$ whose spherical coordinate is given by $\theta_0:=(\pi/2,\ldots,\pi/2,\pi)$. Note that, by \eqref{sigmax0}, we have the following standard deviation function under spherical coordinates,
\bqn{ \label{sigmax}
	\tilde{\sigma}(\theta):=\sigma(x)=\arccos^{\beta} (x_N)
	= \arccos^{\beta} \left(\left(\prod_{i=1}^{N-1}\sin \theta_i\right) \cos \theta_{N}\right),
}
which attains its unique maximum at the interior point $\theta_0\in \Theta$ above. Additionally, it follows from \eqref{eq:SFBM-var} and \eqref{eq:SFBM-cov} that the correlation function of $B_\beta(x)$ becomes
\bqn{\label{corrT}
	r(x,y)= \frac{d^{2\beta}(x,o)+d^{2\beta}(y,o)-d^{2\beta}(x,y)}{2d^\beta(x,o)d^\beta(y,o)}, \quad x, y \in \S^N,
}
whose form under spherical coordinates, denoted by $\tilde{r}(\theta,\varphi)$, can be obtained accordingly.

\begin{remark}
	We choose the starting point of $B_\beta(x)$ at $o=(0,0,\ldots,1,0)\in \R^{N+1}$ to make sure that the maximum of the variance function of $\widetilde{B}_\beta(\theta)$ will be attained at an interior point in $\Theta$. This will simplify a lot the arguments on deriving the asymptotics for the excursion probability. Note that the choice of starting point $o$ does not affect our results since the asymptotics of the excursion probability is determined only by the behavior of the field around the points attaining the maximum of the variance function. \qed
\end{remark}

We first derive a result below showing the local behaviors of the standard deviation and correlation functions of $\widetilde{B}_\beta(\theta)$, the SFBM under spherical coordinates, around $\theta_0$.
\BEL\label{lemA} 
Let $\theta_0=(\pi/2,\ldots,\pi/2,\pi)\in \Theta \subset \mathbb{R}^N$. Then
\begin{equation}\label{Eq:sd sphere}
\begin{split}
\tilde{\sigma}(\theta)= \pi^{\beta}-\beta \pi^{\beta-1}\|\theta-\theta_0\|(1+o(1)), \quad \text{\rm as } \|\theta-\theta_0\|\to 0;
\end{split}
\end{equation}
and
\begin{equation}\label{Eq:corr sphere}
\begin{split}
\tilde{r}(\theta,\varphi)&= 1-\frac{\|\varphi-\theta\|^{2\beta}}{2\pi^{2\beta}}(1+o(1)), \quad \text{\rm as } \|\theta-\theta_0\|\vee\|\varphi-\theta_0\|\to 0.
\end{split}
\end{equation}
\EEL
\begin{proof} \, Note that, as $x_N \downarrow -1$,
\begin{equation}\label{eq:arccos}
\begin{split}
\arccos (x_N)-\arccos(-1)&= \int_{-1}^{x_N} \frac{-1}{\sqrt{1-t^2}}dt=\int_{0}^{1+x_N} \frac{-1}{\sqrt{2t-t^2}}dt\\
&\sim -\int_{0}^{1+x_N} \frac{1}{\sqrt{2t}}dt =-\sqrt{2(1+x_N)}.
\end{split}
\end{equation}
It then follows from \eqref{eq:arccos} and Taylor's expansion that, as $\|\theta-\theta_0\|\to 0$,
\begin{equation*}
\begin{split}
\tilde{\sigma}(\theta)-\pi^{\beta}&=\arccos^{\beta} (\prod_{i=1}^{N-1}\sin \theta_i \cos \theta_{N})-\pi^{\beta}\sim \beta \pi^{\beta-1}\left(\arccos (\prod_{i=1}^{N-1}\sin \theta_i \cos \theta_{N})-\arccos(-1)\right)\\
&\sim -\sqrt{2}\beta \pi^{\beta-1}\sqrt{1+\prod_{i=1}^{N-1}\sin \theta_i \cos \theta_{N}}\sim -\beta \pi^{\beta-1}\sqrt{\sum_{i=1}^{N-1}|\theta_i-\pi/2|^2 + |\theta_N-\pi|^2},
\end{split}
\end{equation*}
yielding \eqref{Eq:sd sphere}.

We derive next the expansion for the correlation function. First note that
\begin{equation}\label{eq:1-corr}
1-r(x,y)=\frac{d^{2\beta}(x,y)-\left(d^{\beta}(x,o)-d^{\beta}(y,o)\right)^2}{2d^\beta(x,o)d^\beta(y,o)}.
\end{equation}
As $x, y \to p=(0,0,\ldots,-1,0)\in \R^{N+1}$, which is equivalent to $\theta, \varphi \to \theta_0$, by Taylor's formula,
\begin{equation*}
\left(d^{\beta}(x,o)-d^{\beta}(y,o)\right)^2\sim \beta^2\pi^{2(\beta-1)}(d(x,o)- d(y,o))^2.
\end{equation*}
Combining this with the triangle inequality such that $|d(x,o)- d(y,o)|\le d(x,y)$, we obtain that for $\beta\in (0,1/2]$,
\[
\left(d^{\beta}(x,o)-d^{\beta}(y,o)\right)^2 = o\left(d^{2\beta}(x,y)\right), \quad \text{as} \ x, y \to p,
\]
implying $1-r(x,y) \sim d^{2\beta}(x,y)/(2\pi^{2\beta})$ by \eqref{eq:1-corr}. Note also that by Lemma 2.1 in Cheng and Xiao \cite{ChengXiaoSphere}, as $x, y \to p$, $d(x,y)\sim \|\theta-\varphi\|$. Therefore,
\[
1-\tilde{r}(\theta,\varphi)=1-r(x,y) \sim \frac{d^{2\beta}(x,y)}{2\pi^{2\beta}} \sim \frac{\|\varphi-\theta\|^{2\beta}}{2\pi^{2\beta}},
\]
yielding \eqref{Eq:corr sphere}.
\end{proof}

For convenience, we present here a simpler version of Theorem 8.2 in Piterbarg \cite{Piterbarg:1996}.
Let $\{X(t), t\in E\}$, where $E\subset \mathbb{R}^N$ is a compact set, be a centered Gaussian random field with variance function  attaining its maximum at the unique point $t_0\in E$. Moreover, there exist non-degenerate $N\times N $ matrices $A$ and $C$, and constants $\eta>0$ and $\alpha\in (0,2]$ such that
\begin{equation}\label{Var}
\sqrt{{\rm Var}(X(t))}=1-\|A(t-t_0)\|^\eta(1+o(1)), \quad \|t-t_0\|\to 0,
\end{equation}
and
\begin{equation}\label{Corr}
{\rm Corr}(X(t), X(s))=1-\|C(t-s)\|^\alpha(1+o(1)), \quad t,s\to t_0.
\end{equation}
Additionally, there exist $\gamma>0$ and $G>0$ such that
\begin{equation}\label{Holder}
\mathbb{E}\left[X(t)-X(s)\right]^2\leq G\|t-s\|^\gamma, \quad s,t\in E.
\end{equation}

\BEL\label{Piterbarg} Let $\{X(t), t\in E\}$, where $E\subset \mathbb{R}^N$ is a compact set, be a centered Gaussian random field with variance function  attaining its maximum at the unique point $t_0\in E$. Assume further that $t_0$ is an inner point of $E$ and (\ref{Var})-(\ref{Holder}) are satisfied. \\
If $\alpha<\eta$, then
\begin{equation}\label{eq:Piterbarg-1}
\pk{\sup_{t\in E}X(t)>u}\sim \mathcal{H}_{\alpha}^N\int_{\mathbb{R}^N} e^{-\|AC^{-1}t\|^\eta}dt\,u^{\frac{2N}{\alpha}-\frac{2N}{\eta}}\Psi(u), \quad u\to\infty.
\end{equation}
If $\alpha=\eta$, then
\begin{equation}\label{eq:Piterbarg-2}
\pk{\sup_{t\in E}X(t)>u}\sim \mathcal{P}_{\alpha}^{\|AC^{-1}t\|^{\alpha}}\Psi(u),\quad u\to\infty.
\end{equation}
If $\alpha>\eta$, then
$$\pk{\sup_{t\in E}X(t)>u}\sim \Psi(u), \quad u\to\infty.$$
\EEL

We are now ready to derive one of our main results as follows.
\BT \label{Thm-Sphere} Let $\{B_\beta (x), x \in \mathbb S^N\}$ be a SFBM, where $\beta\in(0,1/2]$.

(i) If $\beta\in (0,1/2)$, then
\[
	\pk{ \sup_{x\in \mathbb{S}^N} B_\beta(x)>u}
	\sim
	\mathcal{H}_{2\beta}^N\frac{N!\pi^{(2\beta-1/2) N}}{2^{\frac{N}{2\beta}}\beta^{N}\Gamma(N/2 +1)} u^{\frac{(1-2\beta)N}{\beta}}\Psi(\pi^{-\beta}u), \quad u\to\infty,
\]
where $\mathcal{H}_{2\beta}^N$ is the Pickands constant defined in \eqref{eq:Pickands-Const}.

(ii) If $\beta=1/2$, then
\[
	\pk{ \sup_{x\in \mathbb{S}^N} B_\beta(x)>u}
	\sim \mathcal{P}_{{1}}^g \Psi(\pi^{-1/2}u),\quad u\to\infty.
\]
where $\mathcal{P}_{{1}}^g$ is the Piterbarg constant defined in \eqref{eq:Piterbarg-Const} and $g(t)=\sqrt{\sum_{i=1}^Nt_i^2}$, $t\in \R^N$.
\ET

\begin{proof} \, (i) Note that
	\BQNY
	\pk{ \sup_{x\in \mathbb{S}^N} B_\beta(x)>u}=\pk{ \sup_{\theta\in \Theta } \frac{\widetilde{B}_\beta(\theta)}{\pi^\beta}>\frac{u}{\pi^{\beta}}}.
	\EQNY
	It follows from Lemma \ref{lemA} that
	\BQNY
	\begin{split}
		1-\frac{\tilde{\sigma}(\theta)}{\pi^\beta}
		&= \frac{\beta}{\pi} \|\theta-\theta_0\|(1+o(1)), \quad \|\theta-\theta_0\|\to 0,\\
		1-\tilde{r}(\theta,\varphi)&=
		\frac{1}{2\pi^{2\beta}}\|\theta-\varphi\|^{2\beta}(1+o(1)), \quad \|\theta-\theta_0\|\vee\|\varphi-\theta_0\|\to 0.
	\end{split}
	\EQNY
	Applying the identity
	$$\|x-y\|=2\sin\left(\frac{d(x,y)}{2}\right), \quad \forall x,y\in \mathbb{S}^N, $$
	there exists a positive constant $C_1$ such that
	$$d^{2\beta}(x,y)\leq C_1\|x-y\|^{2\beta}, \quad \forall x,y\in \mathbb{S}^N.$$
	Combining this inequality with \eqref{eq:SFBM-var}, there exists a positive constant $C_2$ such that
	$$\E\left[\widetilde{B}_\beta(\theta)-\widetilde{B}_\beta(\varphi)\right]^2= d^{2\beta}(x,y)\leq C_1\|x-y\|^{2\beta}\leq C_2 \|\theta-\varphi\|^{2\beta}, \quad \forall \theta, \varphi\in \Theta.$$
	Therefore, for $\beta\in (0,1/2)$, applying Lemma \ref{Piterbarg} with $\eta=1$, $\alpha=2\beta$, $A=\beta\pi^{-1}I_N$ and $C=2^{-1/(2\beta)}\pi^{-1} I_N$, we obtain
	\BQNY
	\pk{ \sup_{x\in \mathbb{S}^N} B_\beta(x)>u}\sim \mathcal{H}_{2\beta}^N\int_{\mathbb{R}^N}e^{-2^{1/(2\beta)}\beta\|s\|}ds\, v^{\frac{(1-2\beta)N}{\beta}}\Psi(v),
	\EQNY
	where $v=\pi^{-\beta}u$ and $I_N$ is the $N\times N$ identity matrix. Note that $\int_0^\IF r^{N-1}e^{-r}dr=\Gamma(N)$ and
	\begin{equation}\label{eq:area-S}
	\int_{[0,\pi]^{N-2}\times[0,2\pi]} \sin^{N-2}\theta_1\sin^{N-3}\theta_2\dots\sin\theta_{N-2}d\theta_1\cdots d\theta_{N-1} = {\rm Area}(\S^{N-1})=\frac{2\pi^{N/2}}{\Gamma(N/2)},
	\end{equation}
	one can use the spherical coordinate transformation to obtain
	\BQNY\int_{\mathbb{R}^N}e^{-\|s\|}ds=\Gamma(N)\times \frac{2\pi^{N/2}}{\Gamma(N/2)}=\frac{N!\pi^{N/2}}{\Gamma(N/2+1)}.
	\EQNY
	Therefore, as $u\to \infty$,
	\BQNY
	\pk{ \sup_{x\in \mathbb{S}^N} B_\beta(x)>u}\sim\mathcal{H}_{2\beta}^N\frac{N!\pi^{(2\beta-1/2) N}}{2^{\frac{N}{2\beta}}\beta^{N}\Gamma(N/2 +1)} u^{\frac{(1-2\beta)N}{\beta}}\Psi(\pi^{-\beta}u).
	\EQNY
	
	(ii) For $\beta=1/2$, applying again Lemma \ref{Piterbarg}, we have that
	\BQNY
	\pk{ \sup_{x\in \mathbb{S}^N} B_\beta(x)>u}\sim \mathcal{P}_{1}^g\Psi(\pi^{-\beta}u), \quad u\to\infty,
	\EQNY
	where $g(t)=\|t\|$, $t\in \R^N$.
\end{proof}

\section{Excursion Probability on a Geodesic Disc}\label{sec:disc}
In this section, we will study the excursion probability of $B_\beta(x)$ over a geodesic disc on $\S^N$. Without loss of generality, assume that $B_\beta(x)$ starts at $o'=(1,0,\ldots,0)\in \R^{N+1}$ whose spherical coordinate is given by $(0,\ldots,0)\in \Theta \subset \R^N$ according to \eqref{Eq:spherical coordinate}. The standard deviation function of $\widetilde{B}_\beta(\theta)$ now becomes
\bqn{ \label{Eq:sd geodesic ball}
	\tilde{\sigma}(\theta)=\sigma(x)=\arccos^{\beta} \l x, o'\r=\arccos^{\beta} (x_1) =\theta_1^\beta, \quad \theta\in \Theta.}
The geodesic disc on $\S^N$ with radius $a>0$ and center at $o'$ is defined as
\[
T_a=\{x\in\mathbb{S}^N: d(x,o')\le a\}.
\]
Since $d(x,o')=\theta_1$, the set corresponding to $T_a$ under spherical coordinates becomes
$$\Theta_a=[0,a]\times[0, \pi]^{N-2}\times[0, 2 \pi).$$
It is straightforward to check that $\tilde{\sigma}(\theta)$ attains its maximum only at
\[
\{\theta\in \Theta_a: \theta_1=a\}=\{a\}\times[0, \pi]^{N-2}\times[0, 2 \pi),
\]
which is one of the $(N-1)$-dimensional faces of the $N$-dimensional rectangle $\Theta_a$.

\begin{remark}
	We choose the starting point of $B_\beta(x)$ at $o'=(1,0,\ldots,0)\in \R^{N+1}$ to make the variance function of $\widetilde{B}_\beta(\theta)$ have a simple form so that the set attaining the maximum of variance would be easier to handle. Again, the choice of starting point does not affect our results. \qed
\end{remark}

Similarly to Lemma \ref{lemA}, we have the following result describing the local behaviors of the standard deviation and correlation functions of $\widetilde{B}_\beta(\theta)$ around $\Theta_a$.

\BEL\label{lemB} 
Let $\Theta_a=[0,a]\times[0, \pi]^{N-2}\times[0, 2 \pi)$. Then
\begin{equation}\label{eqvar}
\frac{\tilde{\sigma}(\theta)}{a^\beta}=1-\frac{\beta}{a} |a-\theta_1|(1+o(1)), \quad  \theta\in \Theta_a, \ \theta_1\to a;
\end{equation}
and
\begin{equation}\label{eqcorr}
\begin{split}
\tilde{r}(\theta,\varphi)&=1-(1+o(1))\frac{1}{2a^{2\beta}}\bigg[(\varphi_1-\theta_1)^2 + (\sin^2 a)(\varphi_2-\theta_2)^2+\cdots \\
&\quad + \bigg(\sin^2 a\prod_{i=2}^{N-1} \sin^2 \theta_i\bigg)(\varphi_N-\theta_N)^2\bigg]^\beta, \quad \theta,\varphi\in \Theta_a, \ \|\theta-\varphi\|\to 0, \ \theta_1\to a.
\end{split}
\end{equation}
\EEL
\begin{proof}\, Note that \eqref{eqvar} follows immediately from Taylor's formula. By similar arguments in the proof of Lemma \ref{lemA}, we obtain
\[
1-\tilde{r}(\theta,\varphi)=1-r(x,y) \sim \frac{d^{2\beta}(x,y)}{2a^{2\beta}}.
\]
Then \eqref{eqcorr} follows from Lemma 2.1 in Cheng and Xiao \cite{ChengXiaoSphere}.
\end{proof}

Here, we present a result extending both Theorems 7.1 and 8.2 in Piterbarg \cite{Piterbarg:1996}. It is not only useful to prove Theorem \ref{ThA} below, but valuable itself in extreme value theory. The proof is given in the Appendix.

Let $\{X(t), t\in E\}$, where $E=\prod_{i=1}^N[a_i, b_i]$, be a Gaussian random field with continuous trajectories. Its standard deviation function $\sigma_X(t)$ attains the maximum $1$ at the hyperspace $E_0=\{t_1^*\}\times \prod_{i=2}^N[a_i,  b_i]$, where $t_1^*\in [a_1, b_1]$, and satisfies
\BQN\label{var}
\lim_{|t_1-t_1^*|\to 0}\sup_{\substack{t_1\neq t_1^*\\ \hat{t}\in \prod_{i=2}^N[a_i,  b_i]} }\left|\frac{1-\sigma_X(t)}{h(\hat{t})|t_1-t_1^*|^{\gamma}}-1\right|=0,
\EQN
where $\gamma>0$ and $h(\hat{t})$, $\hat{t}\in \prod_{i=2}^N[a_i,  b_i]$, is a positive continuous function with $\hat{t}=(t_2,\dots, t_N)$.
Moreover,
\BQN\label{cor}
\lim_{\delta\to 0, \, u\to\infty}\sup_{\substack{s\neq t, \, s,t\in E_u\\ \|s-t\|\leq \delta} }\left|\frac{1-r(s,t)}{ \left(c_1(t_1-s_1)^2+\sum_{i=2}^Nc_i(\hat{t})(t_i-s_i)^2\right)^\beta}-1\right|=0,
\EQN
where $E_u=\left([t_1^*-((\log u)/u)^{2/\gamma}, t_1^*+((\log u)/u)^{2/\gamma}]\times \prod_{i=2}^N[a_i,  b_i] \right)\cap E$, $\delta>0$, $\beta\in (0,1)$, $c_1>0$ and $c_i(\hat{t})$, $2\leq i\leq N$, are positive and continuous functions over $\prod_{i=2}^N[a_i,  b_i]$.
Additionally, assume that
\BQN\label{cor1}
r(s,t)<1,\quad  s\neq t, \ s,t\in E.
\EQN

\BT\label{extension}
Let $\{X(t), t\in E\}$, where $E=\prod_{i=1}^N[a_i, b_i]$, be a Gaussian random field with continuous trajectories satisfying (\ref{var})-(\ref{cor1}) and let $t_1^*= a_1$ or $b_1$. Then we have, as $u\to\infty$,

(i) for $\beta<\gamma/2$,
\begin{equation}\label{eq:asym-1}
\pk{\sup_{t\in E}X(t)>u}\sim \sqrt{c_1}\Gamma(1/\gamma+1)\mathcal{H}_{2\beta}^N \int_{\hat{t}\in \prod_{i=2}^N[a_i,  b_i] }h^{-1/\gamma}(\hat{t})\prod_{i=2}^N\sqrt{ c_i(\hat{t})}d\hat{t} \, u^{\frac{N}{\beta}-\frac{2}{\gamma}}\Psi(u);
\end{equation}

(ii) for $\beta=\gamma/2$,
\begin{equation}\label{eq:asym-2}
\pk{\sup_{t\in E}X(t)>u}\sim \int_{\hat{t}\in \prod_{i=2}^N[a_i,  b_i] }\widehat{\mathcal{M}}_{2\beta}^{g(t)}\prod_{i=2}^N\sqrt{ c_i(\hat{t})}d\hat{t} \, u^{\frac{N-1}{\beta}}\Psi(u),
\end{equation}
where $g(t)=c_1^{-\beta}h(\hat{t})|t_1|^\gamma$, $t\in \R^N$;

(iii) for $\beta>\gamma/2$,
\begin{equation}\label{eq:asym-3}
\pk{\sup_{t\in E}X(t)>u}\sim \mathcal{H}_{2\beta}^{N-1}\int_{\hat{t}\in \prod_{i=2}^N[a_i,  b_i] }\prod_{i=2}^N\sqrt{ c_i(\hat{t})}d\hat{t} \, u^{\frac{N-1}{\beta}}\Psi(u).
\end{equation}
\ET
\begin{remark}
In Theorem \ref{extension} above, we consider the case when $t_1^*$ is the boundary of the interval $[a_1, b_1]$. If $t_1^*\in (a_1,b_1)$, the following results will be obtained by modifying the proof accordingly. (i) For $\beta<\gamma/2$, replace $\mathcal{H}_{2\beta}^N$ by $2\mathcal{H}_{2\beta}^N$ in the asymptotics in \eqref{eq:asym-1}; (ii) for $\beta=\gamma/2$, replace $\widehat{\mathcal{M}}_{2\beta}^{g(t)}$ by $\mathcal{M}_{2\beta}^{g(t)}$ in the asymptotics in \eqref{eq:asym-2}; (iii) for $\beta>\gamma/2$, the asymptotics in \eqref{eq:asym-2} still holds. \qed
\end{remark}

We formulate our next main result as following.
\BT \label{ThA} Let $\{B_\beta (x), x \in \mathbb S^N\}$ be a SFBM, where $\beta\in(0,1/2]$, and let $T_a=\{x\in\mathbb{S}^N: d(x,o')\le a\}$ with $a\in (0,\pi)$.

(i) If $\beta\in (0,1/2)$, then
\[
	\pk{\sup_{x\in T_a} B_\beta(x) > u }
	\sim \mathcal{H}_{2\beta}^N\frac{N\pi^{N/2}(\sin a)^{N-1}}{2^{\frac{N}{2\beta}}a^{2N-2\beta-1}\beta\Gamma(N/2+1)} u^{\frac{N}{\beta}-2}\Psi(a^{-\beta}u),
	\quad u\to\infty,
\]
where $\mathcal{H}_{2\beta}^N$ is the Pickands constant defined in \eqref{eq:Pickands-Const}.

(ii) If $\beta=1/2$, then
\[ \pk{\sup_{x\in T_a} B_\beta(x) > u }
	\sim \widehat{\mathcal{M}}_{{1}}^{g}\frac{N\pi^{N/2}(\sin a)^{N-1}}{2^{N-1}a^{2(N-1)}\Gamma(N/2+1)} u^{2(N-1)}\Psi(a^{-1/2}u),
	\quad u\to\infty,
\]
where $\widehat{\mathcal{M}}_{{1}}^{g}$ is defined in \eqref{eq:M-Const} and is finite by Lemma \ref{finiteness}, $g(t)=|t_1|$, $t=(t_1, \ldots, t_N)\in \R^N$.
\ET
\begin{proof}\,
 Note that
\BQNY
\pk{ \sup_{x\in T_a} B_\beta(x)>u}=\pk{ \sup_{\theta\in \Theta_a } \frac{\widetilde{B}_\beta(\theta)}{a^\beta}>\frac{u}{a^\beta}},
\EQNY
and we will focus on studying the excursion probability on the right hand side which turns out to be of Euclidean case. It is straightforward that for any $0<\ep<\pi/2$,
\begin{equation}\label{eq:bounds}
\begin{split}
\pk{ \sup_{\theta\in \Theta_a^\ep } \frac{\widetilde{B}_\beta(\theta)}{a^\beta}>\frac{u}{a^\beta}}&\leq \pk{ \sup_{\theta\in \Theta_a } \frac{\widetilde{B}_\beta(\theta)}{a^\beta}>\frac{u}{a^\beta}}\\
&\leq \pk{ \sup_{\theta\in \Theta_a^\ep } \frac{\widetilde{B}_\beta(\theta)}{a^\beta}>\frac{u}{a^\beta}}+\pk{ \sup_{\theta\in \Theta_a\setminus \Theta_a^\ep } \frac{\widetilde{B}_\beta(\theta)}{a^\beta}>\frac{u}{a^\beta}},
\end{split}
\end{equation}
where $$\Theta_a^\ep=[0,a]\times [\ep, \pi-\ep]^{N-2}\times[0, 2 \pi-\ep].$$
Applying Lemma \ref{lemB} and Theorem \ref{extension} with $\beta=\beta$, $\gamma=1$, $h(\hat{\theta})=\beta/a$, $c_1=(2a^{2\beta})^{-1/\beta}$, $c_2(\hat{\theta})=(2a^{2\beta})^{-1/\beta}\sin^2 a$ and $c_j(\hat{\theta})=(2a^{2\beta})^{-1/\beta}(\sin^2 a)\prod_{i=2}^{j-1} \sin^2 \theta_i$ for $3\le j\le N$, we have that for $\beta\in (0,1/2)$,
\begin{equation}\label{eq:asym-ep-1}
\begin{split}
&\pk{ \sup_{\theta\in \Theta_a^\ep } \frac{\widetilde{B}_\beta(\theta)}{a^\beta}>\frac{u}{a^\beta}}\\
&\quad \sim \mathcal{H}_{2\beta}^N \frac{a}{\beta} (2a^{2\beta})^{-N/(2\beta)}(\sin a)^{N-1}a^{2\beta-N}\int_{\hat{\theta}\in \hat{\Theta}_\ep}\prod_{i=2}^{N-1} (\sin\theta_i)^{N-i}d\hat{\theta} \, u^{N/\beta-2}\Psi(\frac{u}{a^\beta});
\end{split}
\end{equation}
and for $\beta=1/2$,
\begin{equation}\label{eq:asym-ep-2}
\begin{split}
&\pk{ \sup_{\theta\in \Theta_a^\ep } \frac{\widetilde{B}_\beta(\theta)}{a^\beta}>\frac{u}{a^\beta}}\\
&\quad \sim \widehat{\mathcal{M}}_{1}^{g} (2a)^{1-N}(\sin a)^{N-1} a^{1-N} \int_{\hat{\theta}\in \hat{\Theta}_\ep}\prod_{i=2}^{N-1} (\sin\theta_i)^{N-i}d\hat{\theta}\, u^{2(N-1)}\Psi(\frac{u}{a^{1/2}}),
\end{split}
\end{equation}
where $g(t)=|t_1|$, $t=(t_1, \ldots, t_N)\in \R^N$, and $$\hat{\Theta}_\ep=[\ep, \pi-\ep]^{N-2}\times[0, 2 \pi-\ep], \quad \hat{\theta}=(\theta_2,\dots, \theta_N).$$

Next we show that the last term in \eqref{eq:bounds} is negligible. Denote by
$$E_0=[0,a-\ep]\times [0, \pi]^{N-2}\times[0, 2 \pi), \quad E_j=[a-\ep, a]\times F_j, \quad 1\leq j\leq n, $$
where  $F_j, 1\leq j\leq n$, is a collection of compact rectangles forming a partition of $ [0, \pi]^{N-2}\times[0, 2 \pi)\setminus \hat{\Theta}_\ep$. Moreover, assume that  $F_j$ and $F_{j'}$ have no common inner point for $j\neq j'$ and the largest edge of $F_j$  has length $L$. Then we have that
\BQNY
\pk{ \sup_{\theta\in \Theta_a\setminus \Theta_a^\ep } \frac{\widetilde{B}_\beta(\theta)}{a^\beta}>\frac{u}{a^\beta}}\leq \sum_{j=0}^n\pk{ \sup_{\theta\in  E_j } \frac{\widetilde{B}_\beta(\theta)}{a^\beta}>\frac{u}{a^\beta}}.
\EQNY
It follows from Lemma \ref{lemB} that there exists $0<\delta<1$ such that $\sup_{\theta\in E_0 } \widetilde{\sigma}(\theta)/a^\beta< 1-\delta$.
By the Borell-TIS inequality \cite{Adler07}, for $u$ sufficiently large,
\BQNY
\pk{ \sup_{\theta\in E_0 } \frac{\widetilde{B}_\beta(\theta)}{a^\beta}>\frac{u}{a^\beta}}\leq \exp\left\{-\frac{(u/a^\beta)^2}{2(1-\delta)^2}\right\}.
\EQNY
By (\ref{eqvar})-(\ref{eqcorr}) and the Slepain inequality, we have that for $\ep>0$ and $L>0$ sufficiently small
\BQNY
\begin{split}
	\pk{ \sup_{\theta\in E_j } \frac{\widetilde{B}_\beta(\theta)}{a^\beta}>\frac{u}{a^\beta}}&\leq \pk{ \sup_{\theta\in E_j } \frac{\widetilde{B}_\beta(\theta)/\tilde{\sigma}(\theta)}{1+\frac{\beta}{2a} |a-\theta_1|}>\frac{u}{a^\beta}}\leq \pk{ \sup_{\theta\in E_j } \frac{Y(C\theta)}{1+\frac{\beta}{2a} |a-\theta_1|}>\frac{u}{a^\beta}},
\end{split}
\EQNY
where $1\leq j\leq n$, $C>2^{-1/(2\beta)}a^{-1}$, and  $Y(t)$ is a centered homogeneous Gaussian random field with continuous trajectories, unit variance and correlation function satisfying
$${\rm Corr}(Y(s),Y(t))=e^{-\|s-t\|^{2\beta}}, \quad s,t\in \R^N.$$
In light of Theorem \ref{extension} with $\beta=\beta$, $\gamma=1$, $h(\hat{\theta})=\beta/(2a)$, $c_1=C^2$, $c_i(\hat{\theta})=C^2,~ 2\leq i\leq n$, we have for $1\leq j\leq n$,
\BQNY
\begin{split}
\pk{ \sup_{\theta\in E_j } \frac{Y(C\theta)}{1+\frac{\beta}{2a} |a-\theta_1|}>\frac{u}{a^\beta}}&\sim\mathcal{H}_{2\beta}^N C^N \frac{2a}{\beta} {\rm mes}(F_j) (\frac{u}{a^\beta})^{\frac{N}{\beta}-2}\Psi(\frac{u}{a^\beta}), \quad \beta<1/2,\\
\pk{ \sup_{\theta\in E_j } \frac{Y(C\theta)}{1+\frac{\beta}{2a} |a-\theta_1|}>\frac{u}{a^\beta}}&\sim\widehat{\mathcal{M}}_{1}^{(4aC)^{-1}|t_1|}   {\rm mes}(F_j) (\frac{u\sqrt{C}}{\sqrt{a}})^{2N-2}\Psi(\frac{u}{a^{1/2}}), \quad \beta=1/2.
\end{split}
\EQNY
Further,
\BQNY
\begin{split}
\sum_{j=1}^n\pk{ \sup_{\theta\in E_j } \frac{Y(C\theta)}{1+\frac{\beta}{2a} |a-\theta_1|}>\frac{u}{a^\beta}}&\sim\mathcal{H}_{2\beta}^N C^N \frac{2a}{\beta} \sum_{j=1}^n{\rm mes}(F_j) (\frac{u}{a^\beta})^{\frac{N}{\beta}-2}\Psi(\frac{u}{a^\beta}), \quad \beta<1/2,\\
\sum_{j=1}^n\pk{ \sup_{\theta\in E_j } \frac{Y(C\theta)}{1+\frac{\beta}{2a} |a-\theta_1|}>\frac{u}{a^\beta}}&\sim\widehat{\mathcal{M}}_{1}^{(4aC)^{-1}|t_1|}  \sum_{j=1}^n {\rm mes}(F_j) (\frac{u\sqrt{C}}{\sqrt{a}})^{2N-2}\Psi(\frac{u}{a^{1/2}}), \quad \beta=1/2.
\end{split}
\EQNY
Note that $\lim_{\ep\to 0} \sum_{j=1}^n {\rm mes}(F_j)=0$ implies that the last term in \eqref{eq:bounds} is negligible.

Applying \eqref{eq:area-S}, one has
$$\lim_{\ep\to 0}\int_{\hat{\theta}\in \hat{\Theta}_\ep}\prod_{i=2}^{N-1} (\sin\theta_i)^{N-i}d\hat{\theta}=\int_{\hat{\theta}\in [0, \pi]^{N-2}\times[0, 2 \pi)}\prod_{i=2}^{N-1} (\sin\theta_i)^{N-i}d\hat{\theta}=\frac{N\pi^{N/2}}{\Gamma(N/2+1)}. $$
Plugging this into \eqref{eq:asym-ep-1} and \eqref{eq:asym-ep-2}, together with \eqref{eq:bounds}, we obtain the desired asymptotic results by letting $\ep\to 0$.
\end{proof}

It is worth mentioning that, when $N=1$, the geodesic disc $T_a$ becomes a circular arc and the maximum of variance is attained at only the two boundary points of $T_a$. Recall Lemma \ref{Piterbarg} and note that if $N=1$ and $t_0$ is a boundary point instead of an interior point, then we can obtain the asymptotics by multiplying the original asymptotics in \eqref{eq:Piterbarg-1} by 1/2 for $\beta\in (0,1/2)$, and by replacing $\mathcal{P}_{\alpha}^{\|AC^{-1}t\|^{\alpha}}$ in \eqref{eq:Piterbarg-2} by $\widehat{\mathcal{M}}_{\alpha}^{\|AC^{-1}t\|^{\alpha}}$ for $\beta=1/2$. Applying these results, together with Lemma \ref{lemB}, similarly to the proof of Theorem \ref{Thm-Sphere}, we have that
\begin{equation}\label{eq:N1}
\begin{split}
\pk{ \sup_{x\in T_a} B_\beta(x)>u}&\sim
\mathcal{H}_{2\beta}^1\frac{a^{2\beta-1}}{2^{\frac{1}{2\beta}-1}\beta} u^{\frac{1}{\beta}-2}\Psi(a^{-\beta}u), \quad \beta\in (0,1/2),\\
\pk{ \sup_{x\in T_a} B_\beta(x)>u}
&\sim 2\widehat{\mathcal{M}}_{{1}}^{g} \Psi(a^{-1/2}u), \quad \beta=1/2,
\end{split}
\end{equation}
where $g(t)=|t|$, $t\in \R$. Then it is easy to check that the asymptotics in \eqref{eq:N1} are exactly the same as those in Theorem \ref{ThA} for $N=1$.

\section*{Acknowledgments} The authors thank Enkelejd Hashorva from University of Lausanne for his helpful discussions and suggestions, and thank anonymous referees for their careful reading and valuable comments.

\section{Appendix}
The following useful lemma can be shown by similar arguments in the proof of Lemma 7.1 in \cite{Piterbarg:1996}. The proof is omitted in this paper.
\BEL\label{lemma} Let $Y(t), t\in \mathbb{R}^N$ be a centered homogeneous Gaussian random field with continuous trajectories, unit variance and correlation function satisfying, with $\beta\in (0,1]$,
$$1-{\rm Corr}(Y(s),Y(t))=\|s-t\|^{2\beta}(1+o(1)), \quad \|s-t\| \to 0.$$
Denote by $\{u_\lambda, \lambda\in \Lambda\}$ a series of function of $u$ with the property that
$$\lim_{u\to\infty}\sup_{\lambda\in \Lambda}\left|\frac{u_\lambda}{u}-1\right|=0.$$
Then for all $b\geq 0$
$$\lim_{u\to\infty}\sup_{\lambda\in \Lambda}\left|\frac{\pk{\sup_{t\in [0, u^{-1/\beta}S]\times \prod_{i=2}^N[a_i,b_i]}\frac{Y(t)}{1+b|t_1|^{2\beta}}>u_\lambda}}{u^{\frac{N-1}{\beta}}\Psi(u_\lambda)}- \mathcal{M}_{2\beta}^{b|t_1|^{2\beta}}([0,S]) \prod_{i=2}^N(b_i-a_i)\right|=0,$$
\EEL
where
$$\mathcal{M}
_{2\beta}^{b|t_1|^{2\beta}}([0,S]) =\lim_{S_1\to\infty}\frac{\mathcal{P}_{2\beta}^{b|t_1|^{2\beta}}([0,S]\times[0,S_1]^{N-1})}{S_1^{N-1}}\in (0,\IF).$$

The following lemma shows the finiteness of the constant $\widehat{\mathcal{M}}_{{1}}^{g}$ in Theorem \ref{ThA}.
\BEL\label{finiteness} For any $\beta\in (0,1]$ and $b>0$,
$$
\widehat{\mathcal{M}}_{2\beta}^{b|t_1|^{2\beta}}:=\lim_{S\to\infty}\mathcal{M}_{2\beta}^{b|t_1|^{2\beta}}([0,S])\in (0,\infty).
$$
\EEL
\begin{proof}\, Let $Y(t)$ be as in Lemma \ref{lemma}. Note that for $0<S_1<\log u$ and $S>0$,
\BQN\label{main1}
A_0(u,S_1)\leq \pk{\sup_{t\in [0, u^{-1/\beta}\log u]\times [0,1]^{N-1}}\frac{Y(t)}{1+b|t_1|^{2\beta}}>u}\leq \sum_{k=0}^{\lfloor (\log u)/S \rfloor+1} A_k(u,S),
\EQN
where
$$A_k(u,S)=\pk{\sup_{t\in [ u^{-1/\beta}kS, u^{-1/\beta}(k+1)S]\times [0,1]^{N-1}}\frac{Y(t)}{1+b|t_1|^{2\beta}}>u}.$$
In light of Lemma \ref{lemma} and subadditivity of $\mathcal{M}_{2\beta}^0([0,S])$, we have
\BQNY
\begin{split}
	A_0(u,S)&\sim \mathcal{M}_{2\beta}^{b|t_1|^{2\beta}}([0,S])u^{\frac{N-1}{\beta}}\Psi(u),\\
	A_k(u,S)&\leq \pk{\sup_{t\in [ u^{-1/\beta}kS, u^{-1/\beta}(k+1)S]\times [0,1]^{N-1}}Y(t)>u(1+u^{-2}b|kS|^{2\beta})}\\
	&\sim \mathcal{M}_{2\beta}^0([0,S])u^{\frac{N-1}{\beta}}\Psi(u(1+u^{-2}b|kS|^{2\beta}))\\
	&\leq CSe^{-b|kS|^{2\beta}}u^{\frac{N-1}{\beta}}\Psi(u), \quad 1\leq k\leq \lfloor (\log u)/S \rfloor+1,
\end{split}
\EQNY
where $C>0$ is a fixed constant. Dividing (\ref{main1}) by $u^{\frac{N-1}{\beta}}\Psi(u)$ on both sides, we have that
$$\mathcal{M}_{2\beta}^{b|t_1|^{2\beta}}([0,S_1])\leq \mathcal{M}_{2\beta}^{b|t_1|^{2\beta}}([0,S]) +\sum_{k=1}^\infty CSe^{-b|kS|^{2\beta}}<\IF.$$
Letting $S_1\to\infty$ leads to
$$\lim_{S_1\to\infty} \mathcal{M}
_{2\beta}^{b|t_1|^{2\beta}}([0,S_1])<\IF,$$
completing the proof.
\end{proof}

\prooftheo{extension} Without loss of generality, we assume that $t_1^*=a_1$, implying $E_u=[a_1,a_1+((\log u)/u)^{2/\gamma}]\times \prod_{i=2}^N[a_i,  b_i]$. Denote by $$F_{1,\ep}=[a_1+\ep, b_1]\times \prod_{i=2}^N[a_i,  b_i], \quad F_{2,\ep}(u)=[a_1+((\log u)/u)^{2/\gamma}, a_1+\ep]\times \prod_{i=2}^N[a_i,  b_i].$$
Then it follows that
\BQN\label{main}
\begin{split}
	\pk{\sup_{t\in E_u}X(t)>u}&\leq \pk{\sup_{t\in E}X(t)>u}\\
	&\leq \pk{\sup_{t\in E_u}X(t)>u}+\pk{\sup_{t\in F_{1,\ep}}X(t)>u}+\pk{\sup_{t\in F_{2,\ep}(u)}X(t)>u}.
\end{split}
\EQN
By (\ref{var}), for $\ep>0$ sufficiently small, there exists a constant $\delta>0$ such that
$\sup_{t\in F_{1,\ep}}\sigma^2(t)<1-\delta.$
By the Borell-TIS inequality \cite{Adler07}, for $u$ large enough,
\BQN\label{upper1}
\pk{\sup_{t\in F_{1,\ep}}X(t)>u}\leq \exp\left\{-\frac{u^2}{2(1-\delta)}\right\}.
\EQN
Moreover, in light of (\ref{cor}), there exists $C>0$ such that for $u$ sufficiently large and $\ep>0$ sufficiently small,
\BQNY
\E\{\left(\overline{X}(t)-\overline{X}(s)\right)\}\leq C\|t-s\|^{2\beta}\leq NC\sum_{i=1}^N|t_i-s_i|^{2\beta}, \quad s,t\in F_{2,\ep}(u),
\EQNY
where $\overline{X}$ is the standardized field of $X$. Additionally, it follows from \eqref{var} that there exists $C_1>0$ such that
$$\sup_{t\in F_{2,\ep}(u)}\sigma^2(t)<1-C_1\left(\frac{\log u}{u}\right)^2.$$
By Theorem 8.1 in \cite{Piterbarg:1996}, we have that, for $u$ sufficiently large,
\BQN\label{upper2}
\pk{\sup_{t\in F_{2,\ep}(u)}X(t)>u}\leq C_2 u^{\frac{N}{\beta}}\Psi\left(\frac{u}{\sqrt{1-C_1\left(\frac{\log u}{u}\right)^2}}\right).
\EQN

We study next $\pk{\sup_{t\in E_u}X(t)>u}$ to derive the exact asymptotics and show that
\[
\pk{\sup_{t\in F_{1,\ep}}X(t)>u} \quad {\rm and} \quad  \pk{\sup_{t\in F_{2,\ep}(u)}X(t)>u}
\]
are negligible as $u\to\infty$. We distinguish three scenarios: $\beta<\gamma/2$, $\beta=\gamma/2$ and $\beta>\gamma/2$.

{\it \textbf{(i) Case $\bm{\beta<\gamma/2}$}}. We first introduce some notation for further analysis. Let
\BQN\label{IM}
I_{k}(u)=[a_1+ku^{-1/\beta}S, a_1+(k+1)u^{-1/\beta}S],\quad  M_u^{\pm}=\left[\frac{((\log u)/u)^{2/\gamma}}{u^{-1/\beta}S}\right]\pm 1.
\EQN
Split $\prod_{i=2}^N[a_i,  b_i]$ into $n^{N-1}$ rectangles with the form $\prod_{i=2}^N\left[a_i+\frac{k_i(b_i-a_i)}{n}, a_i+\frac{(k_i+1)(b_i-a_i)}{n} \right]$ with $n, k_i\in \mathbb{N}$, denoted by  $\{D_j, 1\leq j\leq n^{N-1}\}$. We assume that $D_j$ and $D_{j'}$ have no common inner points for $j\neq j'$. Let
\BQN\label{CJ}
I_{k,j}(u)&=&I_{k}(u)\times D_j, \quad \Lambda^{\pm}=\{(k,j): 0\leq k\leq M_u^{\pm}, 1\leq j\leq n^{N-1}\},\nonumber\\
u_{k,j,\ep}&=&u\left(1+(1-\ep)h_j\inf_{t_1\in I_k(u)}|t_1-t_1^*|^{\gamma}\right),\quad h_j=\inf_{\hat{t}\in D_j }h(\hat{t}),\\ c(j)&=&((c_1+\ep)^{1/2}, (c_{2,j}+\ep)^{1/2}, \dots, (c_{N,j}+\ep)^{1/2}),\quad c_{k,j}=\sup_{\hat{t}\in D_j}c_k(\hat{t}), 2\leq k\leq N.\nonumber
\EQN
Moreover, let $Y(t)$ be a centered homogeneous Gaussian random fields with continuous trajectories, unit variance and correlation function satisfying ${\rm Corr}(Y(s),Y(t))=e^{-\|s-t\|^{2\beta}}$ with $\beta\in (0,1]$.
It follows straightforwardly that
\BQN\label{Bonferroni}
\pi^-(u)-\Sigma(u)\leq \pk{\sup_{t\in E_u}X(t)>u}\leq \pi^+(u),
\EQN
where
$$\pi^\pm(u)=\sum_{(k,j)\in \Lambda^{\pm}}\pk{\sup_{t\in I_{k,j}(u)}X(t)>u},$$
$$\Sigma(u)=\sum_{(k,j)\neq (k',j'), k\leq k',(k,j), (k',j')\in \Lambda^-}\pk{\sup_{t\in I_{k,j}(u)}X(t)>u, \sup_{t\in I_{k',j'}(u)}X(t)>u}.$$
{\it \underline{Asymptotics for $\pi^\pm(u)$}}.
To derive the upper bound, in light of  Slepian inequlaity we have
\BQNY
\begin{split}
	\pk{\sup_{t\in I_{k,j}(u)}X(t)>u}&\leq\pk{\sup_{t\in I_{k,j}(u)}\overline{X}(t)>u_{k,j,\ep}}\leq\pk{\sup_{t\in I_{k,j}(u)}Y(c(j)t)>u_{k,j,\ep}}\\
	&=\pk{\sup_{t\in c(j)I_{0,1}(u)}Y(t)>u_{k,j,\ep}},
\end{split}
\EQNY
where for any $D\subset R^N$, $c(j)D=\{((c_1+\ep)^{1/2}t_1, (c_{2,j}+\ep)^{1/2}t_2, \dots, (c_{N,j}+\ep)^{1/2}t_N): t\in D\}$.
In light of Lemma \ref{lemma}, we have
\BQN\label{pkj}
\begin{split}
&\pk{\sup_{t\in I_{k,j}(u)}X(t)>u}\\
&\quad \leq  \mathcal{H}_{2\beta}^N([0, (c_1+\ep)^{1/2}S])\prod_{i=2}^N(c_{i,j}+\ep)^{1/2}{\rm mes}(D_j)u^{\frac{N-1}{\beta}}\Psi(u_{k,j,\ep})(1+o(1)),
\end{split}
\EQN
as $u\to\infty,$ uniformly with respect to $(k,j)\in\Lambda^+$. Hence, as $u\to\infty$,
\BQNY
\begin{split}
	\sum_{k=0}^{M_u^+}\pk{\sup_{t\in I_{k,j}(u)}X(t)>u}
	\leq  \mathcal{H}_{2\beta}^N([0, (c_1+\ep)^{1/2}S])\prod_{i=2}^N(c_{i,j}+\ep)^{1/2}{\rm mes}(D_j)u^{\frac{N-1}{\beta}}\sum_{k=0}^{M_u^+}\Psi(u_{k,j,\ep}).
\end{split}
\EQNY
Noting that, as $u\to\infty$,
\BQNY
\begin{split}
\sum_{k=0}^{M_u^+}\Psi(u_{k,j,\ep})&\leq \Psi(u)\sum_{k=0}^{M_u^+}
	e^{-(1-\ep)h_ju^2|ku^{-1/\beta}S|^{\gamma}}\\
&\leq \Psi(u)
	\left((1-\ep)^{1/\gamma}h_j^{1/\gamma}u^{2/\gamma-1/\beta}S\right)^{-1}\\
	&\quad \times\sum_{k=0}^{M_u^+}
	e^{-|k(1-\ep)^{1/\gamma}h_j^{1/\gamma}u^{2/\gamma-1/\beta}S|^{\gamma}}\times (1-\ep)^{1/\gamma}h_j^{1/\gamma}u^{2/\gamma-1/\beta}S\\
&\leq \Psi(u)
	\left((1-\ep)^{1/\gamma}h_j^{1/\gamma}u^{2/\gamma-1/\beta}S\right)^{-1}\int_0^\infty e^{-|t|^\gamma} dt,
\end{split}
\EQNY
we have
\BQNY
\begin{split}
	&\quad \sum_{k=0}^{M_u^+}\pk{\sup_{t\in I_{k,j}(u)}X(t)>u}\\
	&\leq \frac{\mathcal{H}_{2\beta}^N([0, (c_1+\ep)^{1/2}S])}{S}\Gamma(1/\gamma+1)(1-\ep)^{-1/\gamma}h_j^{-1/\gamma}\prod_{i=2}^N(c_{i,j}+\ep)^{1/2}{\rm mes}(D_j)
	u^{\frac{N}{\beta}-\frac{2}{\gamma}}\Psi(u)\\
	&\leq \mathcal{H}_{2\beta}^N\Gamma(1/\gamma+1)(c_1+\ep)^{1/2}(1-\ep)^{-1/\gamma}h_j^{-1/\gamma}\prod_{i=2}^N(c_{i,j}+\ep)^{1/2}{\rm mes}(D_j)
	u^{\frac{N}{\beta}-\frac{2}{\gamma}}\Psi(u),
\end{split}
\EQNY

as $u\to\infty$ and $S\to\infty$.
Furthermore,
\begin{equation}\label{pi+}
\begin{split}
\pi^+(u)&=\sum_{j=1}^{n^{N-1}}\sum_{k=0}^{M_u^+}\pk{\sup_{t\in I_{k,j}(u)}X(t)>u}\\
&\leq  \mathcal{H}_{2\beta}^N\frac{\Gamma(1/\gamma+1)(c_1+\ep)^{1/2}}{(1-\ep)^{1/\gamma}}u^{\frac{N}{\beta}-\frac{2}{\gamma}}\Psi(u)\sum_{j=0}^{n^{N-1}}\left(h_j^{-1/\gamma}\prod_{i=2}^N(c_{i,j}+\ep)^{1/2}{\rm mes}(D_j)\right)\\
&\sim  \mathcal{H}_{2\beta}^N \Gamma(1/\gamma+1)c_1^{1/2}u^{\frac{N}{\beta}-\frac{2}{\gamma}}\Psi(u)\int_{\hat{t}\in \prod_{i=2}^N [a_i,b_i] }h^{-1/\gamma}(\hat{t})\prod_{i=2}^N c_i^{1/2}(\hat{t})d\hat{t},
\end{split}
\end{equation}
as $u\to\infty, n\to\infty, \ep\to 0$. Analogously,
we can show that
\BQN\label{pi-}
\pi^-(u)\sim \mathcal{H}_{2\beta}^N \Gamma(1/\gamma+1)c_1^{1/2}u^{\frac{N}{\beta}-\frac{2}{\gamma}}\Psi(u)\int_{\hat{t}\in \prod_{i=2}^N [a_i,b_i] }h^{-1/\gamma}(\hat{t})\prod_{i=2}^N c_i^{1/2}(\hat{t})d\hat{t},
\EQN
as $u\to\infty, n\to\infty$. Next we show that $\Sigma(u)$ is negligible compared with $\pi^-(u)$ as $u\to\infty$. For this, denote by
\BQNY
\begin{split}
	\Lambda^-_1&=\{(k,j, k',j'): (k,j), (k',j')\in \Lambda^-, k\leq k', D_j\cap D_{j'}=\emptyset\},\\
	\Lambda^-_2&=\{(k,j, k',j'):(k,j), (k',j')\in \Lambda^-,  k\leq k'\leq k+1, D_j\cap D_{j'}\neq \emptyset, j\neq j'\},\\
	\Lambda^-_3&=\{(k,j, k',j'): (k,j), (k',j')\in \Lambda^-,   k+1<k',  D_j\cap D_{j'}\neq \emptyset\},\\
	\Lambda^-_4&=\{(k,j, k',j'): (k,j), (k',j')\in \Lambda^-,  k'=k+1, j=j' \}.
\end{split}
\EQNY
Then it follows that $\Sigma(u)\leq \sum_{i=1}^4\Sigma_i(u)$, where
$$\Sigma_i(u)=\sum_{(k,j,k',j')\in \Lambda^-_i}\pk{\sup_{t\in I_{k,j}(u)}X(t)>u, \sup_{t\in I_{k',j'}(u)}X(t)>u}.$$
{\it\underline{Upper bound for $\Sigma_1(u)$}}. Note that
\BQNY
\Sigma_1(u)\leq\sum_{(k,j,k',j')\in \Lambda^-_1}\pk{\sup_{s\in I_{k,j}(u), t\in I_{k,j}(u)}X(s)+X(t)>2u},
\EQNY
and by (\ref{var}) and (\ref{cor1}), there exists $0<\delta<1$ such that
\BQNY
{\rm Var}\left(X(s)+X(t)\right)= \sigma^2(s)+\sigma^2(t)+2\sigma(s)2\sigma(t)r(s,t)<4-\delta.
\EQNY
It follows from the Borell-TIS inequality \cite{Adler07} that, as $u\to\infty$,
\BQN\label{sigma1}
\Sigma_1(u)\leq\sum_{(k,j,k',j')\in \Lambda^-_1} e^{-\frac{\left(2u-\E\left(\sup_{t\in E}X(t)\right)\right)^2}{2(4-\delta)}}\leq (n^{N-1} M_u^+)^2e^{-\frac{\left(2u-\E\left(\sup_{t\in E}X(t)\right)\right)^2}{2(4-\delta)}}=o(\pi^-(u)).
\EQN
{\it\underline{Upper bound for $\Sigma_2(u)$}}.   For $(k,j,k',j')\in \Lambda^-_2$, without loss of generality, we assume that
\BQNY
\begin{split}
	D_j&=\prod_{i=2}^N\left[a_i+\frac{k_i(b_i-a_i)}{n}, a_i+\frac{(k_i+1)(b_i-a_i)}{n} \right],\\
	D_{j'}&=\left[a_2+\frac{(k_2+1)(b_2-a_2)}{n}, a_2+\frac{(k_2+2)(b_2-a_2)}{n} \right]\\
	&\quad \times\prod_{i=3}^N\left[a_i+\frac{k_i(b_i-a_i)}{n}, a_i+\frac{(k_i+1)(b_i-a_i)}{n} \right].
\end{split}
\EQNY
Split $D_{j'}$ into two parts:
\BQNY
\begin{split}
	D_{j'}^{(1)}&=\left[a_2+\frac{(k_2+1)(b_2-a_2)}{n}, a_2+\frac{(k_2+1)(b_2-a_2)}{n}+\frac{b_2-a_2}{n^{2}} \right]\\
	&\quad \times\prod_{i=3}^N\left[a_i+\frac{k_i(b_i-a_i)}{n}, a_i+\frac{(k_i+1)(b_i-a_i)}{n} \right],\\
	D_{j'}^{(2)}&=\left[a_2+\frac{(k_2+1)(b_2-a_2)}{n}+\frac{b_2-a_2}{n^{2}} , a_2+\frac{(k_2+2)(b_2-a_2)}{n} \right]\\
	&\quad \times\prod_{i=3}^N\left[a_i+\frac{k_i(b_i-a_i)}{n}, a_i+\frac{(k_i+1)(b_i-a_i)}{n} \right].
\end{split}
\EQNY
Then it follows that
\BQNY
\begin{split}
&\pk{\sup_{t\in I_{k,j}(u)}X(t)>u, \sup_{t\in I_{k',j'}(u)}X(t)>u}\\
&\quad \leq \pk{ \sup_{t\in I_{k',j'}^{(1)}(u)}X(t)>u}+\pk{\sup_{t\in I_{k,j}(u)}X(t)>u, \sup_{t\in I_{k',j'}^{(2)}(u)}X(t)>u}
\end{split}
\EQNY
with $I_{k'j'}^{(l)}(u)=I_{k'}(u)\times D_{j'}^{(l)}, l=1,2$.
By Lemma \ref{lemma} and (\ref{pkj}), we have as $u\to\infty$,
\BQNY
\pk{ \sup_{t\in I_{k',j'}^{(1)}(u)}X(t)>u}\leq C\frac{{\rm mes}(D_{j'}^{(1)})}{{\rm mes}(D_{j'})}\pk{ \sup_{t\in I_{k',j'}(u)}X(t)>u},
\EQNY
where $C>0$ is a constant independent of $k'$ and $j'$.
Using the fact that $D_{j'}$ has at most $3^{N-1}$ neighbors and
$$\lim_{n\to\infty}\sup_{1\leq j'\leq n^{N-1}}\frac{{\rm mes}(D_{j'}^{(1)})}{{\rm mes}(D_{j'})}=0,$$
we have
\BQNY
\begin{split}
	\sum_{(k,j,k',j')\in \Lambda^-_2}\pk{ \sup_{t\in I_{k',j'}^{(1)}(u)}X(t)>u}&\leq C\sum_{(k,j,k',j')\in \Lambda^-_2}\frac{{\rm mes}(D_{j'}^{(1)})}{{\rm mes}(D_{j'})}\pk{ \sup_{t\in I_{k',j'}(u)}X(t)>u}\\
	&\leq 3^{N-1}C\sum_{(k',j')\in \Lambda^-}\frac{{\rm mes}(D_{j'}^{(1)})}{{\rm mes}(D_{j'})}\pk{ \sup_{t\in I_{k',j'}(u)}X(t)>u}\\
	&=o(\pi^-(u)),\quad u\to\infty, n\to\infty.
\end{split}
\EQNY
Using the same argument as in (\ref{sigma1}), we have
\BQNY
\sum_{(k,j,k',j')\in \Lambda^-_2}\pk{\sup_{t\in I_{k,j}(u)}X(t)>u,  \sup_{t\in I_{k',j'}^{(2)}(u)}X(t)>u}=o(\pi^-(u)), \quad u\to\infty.
\EQNY
Hence,
$$\Sigma_2(u)=o(\pi^-(u)), \quad u\to\infty.$$
{\it \underline{Upper bound for $\Sigma_3(u)$}}.
Let
\begin{equation*}
\begin{split}
J_{l}(u)&=\prod_{i=2}^N[l_iu^{-1/\beta}S, (l_i+1)u^{-1/\beta}S], \quad \text{with}\quad l=(l_2,\dots,l_N), \quad \Xi_j=\{l: D_j\bigcap J_{l}(u)\neq \emptyset\},\\
J_{k,l}(u)&=I_k(u)\times J_l(u).
\end{split}
\end{equation*}
Then
\BQNY
\begin{split}
	&\pk{\sup_{t\in I_{k,j}(u)}X(t)>u, \sup_{t\in I_{k',j'}(u)}X(t)>u}\\
	&\quad \leq \pk{\sup_{t\in I_{k,j}(u)}\overline{X}(t)>u_{k,j,\ep}, \sup_{t\in I_{k',j'}(u)}\overline{X}(t)>u_{k',j',\ep}}\\
	&\quad \leq \sum_{l\in \Xi_j, l'\in \Xi_{j'}}\pk{\sup_{t\in J_{k,l}(u)}\overline{X}(t)>u_{k,j,\ep}, \sup_{t\in J_{k',l'}(u)}\overline{X}(t)>u_{k',j',\ep}}\\
	&\quad =   \sum_{l\in \Xi_j, l'\in \Xi_{j'}}\pk{\sup_{t\in J_{k,l}(1)}\overline{X}(u^{-1/\beta}t)>u_{k,j,\ep}, \sup_{t\in J_{k',l'}(1)}\overline{X}(u^{-1/\beta}t)>u_{k',j',\ep}}.
\end{split}
\EQNY
In view of \eqref{cor}, there exist $C_1, C_2>0$ such that for $u$ and $n$ sufficiently large
\BQNY
\begin{split}
	C_1\sum_{i=1}^N|s_i-t_i|^{2\beta}\leq u^2(1-{\rm Corr}(\overline{X}(u^{-1/\beta}s),\overline{X}(u^{-1/\beta}t)))\leq C_2\sum_{i=1}^N|s_i-t_i|^{2\beta},
\end{split}
\EQNY
and
$${\rm Corr}(\overline{X}(u^{-1/\beta}s),\overline{X}(u^{-1/\beta}t))\geq 1/2,$$
hold for all $s,t\in [a_1u^{1/\beta},a_1u^{1/\beta}+((\log u)/u)^{2/\gamma}u^{1/\beta}]\times \bigcup_{D_j\cap D_{j'}\neq \emptyset}u^{1/\beta}D_{j'}, 1\leq j\leq n^{N-1}$ with $u^{1/\beta}D_{j'}=\{u^{1/\beta}\hat{t}: \hat{t}\in D_{j'}  \}$.
Thus in light of Corollary 3.1 in \cite{Debicki:2017}, there exist $\mathcal{C}, \mathcal{C}_1>0$ such that for $u$ and $n$ sufficiently large, $l\in \Xi_j, l'\in \Xi_{j'}$, $D_j\cap D_{j'}\neq \emptyset$, $0\leq k,k'\leq M_u^+$ and $|k'-k-1|\geq 1$,
\BQN\label{new}
&\pk{\sup_{t\in J_{k,l}(1)}\overline{X}(u^{-1/\beta}t)>u_{k,j,\ep}, \sup_{t\in J_{k',l'}(1)}\overline{X}(u^{-1/\beta}t)>u_{k',j',\ep}}\nonumber\\
&\quad \leq \mathcal{C} S^{2N}e^{-\mathcal{C}_1S^{2\beta}(|k'-k-1|^{2\beta}+\|l-l'\|^{2\beta})}\Psi(u_{k,k',j,j'\ep}),
\EQN
where $$u_{k,k',j,j'\ep}=\min(u_{k,j,\ep}, u_{k',j',\ep}).$$
We have
\BQNY
\begin{split}
	\Sigma_3(u)&\leq \sum_{(k,j,k',j')\in \Lambda^-_3}\pk{\sup_{t\in I_{k,j}(u)}X(t)>u, \sup_{t\in I_{k',j'}(u)}X(t)>u}\\
	&\leq \sum_{ (k,j)\in \Lambda^-, |k-k'-1|\geq 1, D_j\bigcap D_{j'}\neq \emptyset}\pk{\sup_{t\in I_{k,j}(u)}X(t)>u, \sup_{t\in I_{k',j'}(u)}X(t)>u}\\
&\leq\sum_{ (k,j)\in \Lambda^-, |k-k'-1|\geq 1, D_j\bigcap D_{j'}\neq \emptyset}\sum_{l\in \Xi_j, l'\in \Xi_{j'}}\\
	&\quad \pk{\sup_{t\in J_{k,l}(1)}\overline{X}(u^{-1/\beta}t)>u_{k,j,\ep}, \sup_{t\in J_{k',l'}(1)}\overline{X}(u^{-1/\beta}t)>u_{k',j',\ep}}.
\end{split}
\EQNY
For  $(k,j)\in \Lambda^-$, it follows from (\ref{new}) that
\BQNY
\begin{split}
 &\sum_{|k-k'-1|\geq 1, D_j\bigcap D_{j'}\neq \emptyset}\sum_{l\in \Xi_j, l'\in \Xi_{j'}}
	 \pk{\sup_{t\in J_{k,l}(1)}\overline{X}(u^{-1/\beta}t)>u_{k,j,\ep}, \sup_{t\in J_{k',l'}(1)}\overline{X}(u^{-1/\beta}t)>u_{k',j',\ep}}\\
	&\leq \sum_{ |k-k'-1|\geq 1, D_j\bigcap D_{j'}\neq \emptyset}\sum_{l\in \Xi_j, l'\in \Xi_{j'}}\mathcal{C} S^{2N}e^{-\mathcal{C}_1S^{2\beta}(|k'-k-1|^{2\beta}+\|l-l'\|^{2\beta})}\Psi(u_{k,k',j,j',\ep})\\
	&\leq \sum_{l\in \Xi_j}C_3 S^{2N}e^{-C_4S^{2\beta}}\Psi(u_{k,j,\ep})\\
	&\leq C_3 S^{2N-1}u^{\frac{N-1}{\beta}}e^{-C_4S^{2\beta}}\Psi(u_{k,j,\ep}),
	\quad u\to\infty,
\end{split}
\EQNY
where $C_3$ and $C_4$ are two positive constants.
Hence
\BQNY
\begin{split}
\Sigma_3(u)&\leq \sum_{ (k,j)\in \Lambda^-}C_3 S^{2N-1}u^{\frac{N-1}{\beta}}e^{-C_4S^{2\beta}}\Psi(u_{k,j,\ep})\\
	&= C_3 S^{2N-1}e^{-C_4S^{2\beta}} \sum_{ (k,j)\in \Lambda^-}u^{\frac{N-1}{\beta}}\Psi(u_{k,j,\ep})=o(\pi^-(u)), \quad u\to\infty, S\to\infty.
\end{split}
\EQNY
{\it \underline{Upper bound for $\Sigma_4(u)$}}.
Observe that
\BQNY
\begin{split}
	&\pk{\sup_{t\in I_{k,j}(u)}X(t)>u, \sup_{t\in I_{k+1,j}(u)}X(t)>u}\\
	&\quad=\pk{\sup_{t\in I_{k,j}(u)\bigcup I_{k+1,j}(u)}X(t)>u}-\pk{\sup_{t\in I_{k,j}(u)}X(t)>u}-\pk{\sup_{t\in I_{k+1,j}(u)}X(t)>u}.
\end{split}
\EQNY
Thus in light of  (\ref{pi+}) and (\ref{pi-}), it follows that
\BQNY
\begin{split}
	\Sigma_4(u)&\leq \sum_{(k,j)\in \Lambda^-}\pk{\sup_{t\in I_{k,j}(u)}X(t)>u, \sup_{t\in I_{k+1,j}(u)}X(t)>u}\\
	&\leq \sum_{(k,j)\in \Lambda^-} \pk{\sup_{t\in I_{k,j}(u)\bigcup I_{k+1,j}(u)}X(t)>u}-\sum_{(k,j)\in \Lambda^-} \pk{\sup_{t\in I_{k,j}(u)}X(t)>u}\\
	&\quad -\sum_{(k,j)\in \Lambda^-} \pk{\sup_{t\in I_{k+1,j}(u)}X(t)>u}\\
	&=2\pi^+(u)(1+o(1))-2\pi^-(u)(1+o(1))=o(\pi^-(u)), \quad u\to\infty, S\to\infty.
\end{split}
\EQNY
Therefore we conclude that
\BQNY
\Sigma(u)=o(\pi^-(u)), \quad u\to\infty, S\to\infty,
\EQNY
together with (\ref{Bonferroni}), (\ref{pi+}) and (\ref{pi-}), yielding that as $u\to\infty$,
$$\pk{\sup_{t\in E_u}X(t)>u}\sim \mathcal{H}_{2\beta}^N \Gamma(1/\gamma+1)c_1^{1/2}u^{\frac{N}{\beta}-\frac{2}{\gamma}}\Psi(u)\int_{\hat{t}\in \prod_{i=2}^N [a_i,b_i] }h^{-1/\gamma}(\hat{t})\prod_{i=2}^N c_i^{1/2}(\hat{t})d\hat{t}.$$
Inserting the above asymptotics, (\ref{upper1}) and (\ref{upper2}) into (\ref{main}) establishes the claim.

{\it \textbf{(ii) Case $\bm{\beta=\gamma/2}$}}. It follows that
\BQN\label{Bonferroni1}
\pi_1^-(u)-\Sigma_5(u)\leq \pk{\sup_{t\in E_u}X(t)>u}\leq \pi_1^-(u)+ \pi_1^+(u),
\EQN
where
$$\pi_1^-(u)=\sum_{j=1}^{n^{N-1}}\pk{\sup_{t\in I_{0,j}(u)}X(t)>u},\quad \pi_1^+(u)=\sum_{k=1}^{M_u^+}\sum_{j=1}^{n^{N-1}}\pk{\sup_{t\in I_{k,j}(u)}X(t)>u},$$
$$\Sigma_5(u)=\sum_{1\leq j<j'\leq n^{N-1}}\pk{\sup_{t\in I_{0,j}(u)}X(t)>u, \sup_{t\in I_{0,j'}(u)}X(t)>u},$$
with $I_{k,j}$ and $M_u^{\pm}$ being defined in (\ref{IM}) and (\ref{CJ}).\\
{\it \underline{Asymptotics of $\pi_1^-(u)$}}. By (\ref{var}) and Slepain inequality,
\BQNY
\begin{split}
	\pi_1^-(u)&\leq \sum_{j=1}^{n^{N-1}}\pk{\sup_{t\in I_{0,j}(u)}\frac{\overline{X}(t)}{1+h_j|t_1-t_1^*|^{\gamma}}>u}\leq \sum_{j=1}^{n^{N-1}}\pk{\sup_{t\in I_{0,j}(u)}\frac{Y(c(j)t)}{1+h_j|t_1-t_1^*|^{\gamma}}>u}\\
	&\leq \sum_{j=1}^{n^{N-1}}\pk{\sup_{t\in c(j)I_{0,j}(u)}\frac{Y(t)}{1+h_j(c_1+\ep)^{-\gamma/2}|t_1-t_1^*|^{\gamma}}>u},
\end{split}
\EQNY
where $h_j$ and $c(j)$ are given in (\ref{CJ}).
In light of Lemma \ref{lemma} and Lemma \ref{finiteness}, we have that
\BQN\label{pi1-}
\begin{split}
	\pi_1^-(u)&\leq \sum_{j=1}^{n^{N-1}} \mathcal{M}_{2\beta}^{h_j(c_1+\ep)^{-\gamma/2}|t_1|^\gamma}[0,(c_1+\ep)^{1/2}S]\prod_{i=2}^N(c_{i,j}+\ep)^{1/2}{\rm mes}(D_j)u^{\frac{N-1}{\beta}}\Psi(u)\\
	&\sim \int_{\hat{t}\in \prod_{i=2}^N[a_i,  b_i] }\widehat{\mathcal{M}}_{2\beta}^{c_1^{-\beta}h(\hat{t})|t_1|^\gamma}\prod_{i=2}^N\sqrt{ c_i(\hat{t})}d\hat{t} u^{\frac{N-1}{\beta}}\Psi(u),
\end{split}
\EQN
as $u\to\infty, n\to\infty, \ep\to 0$.
Similarly, we can show that
\BQNY
\pi_1^-(u)\geq \int_{\hat{t}\in \prod_{i=2}^N[a_i,  b_i] }\widehat{\mathcal{M}}_{2\beta}^{c_1^{-\beta}h(\hat{t})|t_1|^\gamma}\prod_{i=2}^N\sqrt{ c_i(\hat{t})}d\hat{t} u^{\frac{N-1}{\beta}}\Psi(u)(1+o(1)), \quad u\to\infty, S\to\infty.
\EQNY
Hence
\BQNY
\pi_1^-(u)\sim\int_{\hat{t}\in \prod_{i=2}^N[a_i,  b_i] }\widehat{\mathcal{M}}_{2\beta}^{c_1^{-\beta}h(\hat{t})|t_1|^\gamma}\prod_{i=2}^N\sqrt{ c_i(\hat{t})}d\hat{t} u^{\frac{N-1}{\beta}}\Psi(u),\quad u\to\infty, S\to\infty.
\EQNY
{\it \underline{Upper bound for $\pi_1^+(u)$}}. In view of (\ref{pkj}), we have, as $u\to\infty$,
\BQNY
\begin{split}
	\pi_1^+(u)&\leq \sum_{k=1}^{M_u^+}\sum_{j=1}^{n^{N-1}}\pk{\sup_{t\in I_{k,j}(u)}\overline{X}(t)>u_{k,j,\ep}}\\
	&\leq \sum_{k=1}^{M_u^+}\sum_{j=1}^{n^{N-1}}\mathcal{H}_{2\beta}^N([0, (c_1+\ep)^{1/2}S])\prod_{i=2}^N(c_{i,j}+\ep)^{1/2}{\rm mes}(D_j)u^{\frac{N-1}{\beta}}\Psi(u_{k,j,\ep}).
\end{split}
\EQNY
Note that as $u\to\infty$,
\BQNY
\begin{split}
	\sum_{k=1}^{M_u^+}\Psi(u_{k,j,\ep})&\sim \Psi(u)\sum_{k=1}^{M_u^+}e^{-(1-\ep)h_j|kS|^{\gamma}}\leq \Psi(u)e^{-CS^{\gamma}},
\end{split}
\EQNY
where $C>0$ is a positive constant. It follows that
\BQN\label{pi1+}
\begin{split}
	\pi_1^+(u)
	&\leq Se^{-CS^{\gamma}}\sum_{j=1}^{n^{N-1}}\frac{\mathcal{H}_{2\beta}^N([0, (c_1+\ep)^{1/2}S])}{S}\prod_{i=2}^N(c_{i,j}+\ep)^{1/2}{\rm mes}(D_j)u^{\frac{N-1}{\beta}}\Psi(u)\\
	&=o(\pi_1^-(u)), \quad u\to\infty, S\to\infty.
\end{split}
\EQN
{\it \underline{Upper bound for $\Sigma_5(u)$}}.
Let
\begin{equation*}
\begin{split}
\Lambda_5&=\{(j,j'): 1\leq j<j'\leq n^{N-1}, D_j\bigcap D_{j'}=\emptyset\},\\
\Lambda_6&=\{(j,j'): 1\leq j<j'\leq n^{N-1}, D_j\bigcap D_{j'}\neq\emptyset\}.
\end{split}
\end{equation*}
Then
\BQNY
\Sigma_5(u)\leq \Sigma_6(u)+\Sigma_7(u),
\EQNY
with
$$\Sigma_i(u)=\sum_{(j,j')\in \Lambda_{i-1}}\pk{\sup_{t\in I_{0,j}(u)}X(t)>u, \sup_{t\in I_{0,j'}(u)}X(t)>u}, \quad i=6,7.$$
Using same arguments as in those to get the upper bounds of $\Sigma_1(u)$ and $\Sigma_2(u)$, we can show that  $\Sigma_i(u)=o(\pi_1^-(u)), i=6,7$, as $u\to\infty$ and $n\to\infty$.
Hence
$$\pk{\sup_{t\in E_u}X(t)>u}\sim \int_{\hat{t}\in \prod_{i=2}^N[a_i,  b_i] }\widehat{\mathcal{M}}_{2\beta}^{c_1^{-\beta}h(\hat{t})|t_1|}\prod_{i=2}^N\sqrt{ c_i(\hat{t})}d\hat{t} u^{\frac{N-1}{\beta}}\Psi(u), \quad u\to\infty,$$
together with (\ref{main})-(\ref{upper2}), establishing the claim.

{\it \textbf{(iii) Case $\bm{\beta>\gamma/2}$}}. Observe that
\BQNY
\pi_2(u)\leq \pk{\sup_{t\in E_u}X(t)>u}\leq \pi_1^-(u)+\pi_1^+(u),
\EQNY
with $\pi_1^-(u)$ and $\pi_1^+(u)$ defined in (\ref{Bonferroni1}) and
$$\pi_2(u)=\pk{\sup_{\hat{t}\in \prod_{i=2}^N [a_i,b_i] }X(t_1^*,\hat{t})>u}.$$
\textit{\underline{Upper bound of $\pi_1^\pm(u)$}}. By (\ref{var}) and (\ref{pi1-}) , it follows that for any $q, \ep>0$  and sufficiently large $S_1$, as $u\to\infty, n\to\infty$,
	\BQNY
	\begin{split}
		\pi_1^-(u)&\leq \sum_{j=1}^{n^{N-1}}\pk{\sup_{t\in I_{0,j}(u)}\frac{\overline{X}(t)}{1+q|t_1-t_1^*|^{2\beta}}>u}\\
		&\sim \mathcal{M}_{2\beta}^{c_1^{-\beta}q|t_1|^{2\beta}}[0,S]\int_{\hat{t}\in \prod_{i=2}^N[a_i,  b_i] }\prod_{i=2}^N\sqrt{ c_i(\hat{t})}d\hat{t} u^{\frac{N-1}{\beta}}\Psi(u)\\
		&\leq (1+\ep)\frac{\mathcal{P}_{2\beta}^{c_1^{-\beta}q|t_1|^{2\beta}}([0,S]\times[0,S_1]^{N-1})}{S_1^{N-1}}\int_{\hat{t}\in \prod_{i=2}^N[a_i,  b_i] }\prod_{i=2}^N\sqrt{ c_i(\hat{t})}d\hat{t} u^{\frac{N-1}{\beta}}\Psi(u).
	\end{split}
	\EQNY
	By the fact that
	$$\lim_{S_1\to\infty}\lim_{q\to\infty}\frac{\mathcal{P}_{2\beta}^{c_1^{-\beta}q|t_1|^{2\beta}}([0,S]\times[0,S_1]^{N-1})}{S_1^{N-1}}=
	\lim_{S_1\to\infty}\frac{\mathcal{H}_{2\beta}^N(\{0\}\times[0,S_1]^{N-1})}{S_1^{N-1}}=\mathcal{H}_{2\beta}^{N-1},$$
	we have that as $u\to\infty$,  $n\to\infty$, $q\to\infty$, $S_1\to\infty$ and $\ep\to 0$,
	\BQNY
	\begin{split}
		\pi_1^-(u)
		\leq \mathcal{H}_{2\beta}^{N-1}\int_{\hat{t}\in \prod_{i=2}^N[a_i,  b_i] }\prod_{i=2}^N\sqrt{ c_i(\hat{t})}d\hat{t} u^{\frac{N-1}{\beta}}\Psi(u).
	\end{split}
	\EQNY
	Using the same argument as in  (\ref{pi1+}), we have that
	$$\pi_1^+(u)=o(\pi_1^-), \quad u\to\infty, S\to\infty.$$
	\textit{\underline{Asymptotics of $\pi_2(u)$}}. Note that $X(t_1^*,\hat{t})$ is a Gaussian random field with unit variance and correlation function satisfying
	\BQNY
	\lim_{\delta\to 0}\sup_{\hat{s}\neq \hat{t}, \hat{s},\hat{t}\in \prod_{i=2}^N[a_i,  b_i], |\hat{t}-\hat{s}|\leq \delta }\left|\frac{1-{\rm Corr}(X(t_1^*,\hat{t}), X(t_1^*,\hat{s}))}{ \left(\sum_{i=2}^Nc_i(\hat{t})(t_i-s_i)^2\right)^\beta}-1\right|=0,
	\EQNY
	and $${\rm Corr}(X(t_1^*,\hat{t}), X(t_1^*,\hat{s}))<1, \quad \hat{t}\neq \hat{s}, \hat{t}, \hat{s}\in  \prod_{i=2}^N[a_i,  b_i].$$
	By Theorem 7.1 in \cite{Piterbarg:1996}, we have
	$$\pi_2(u)\sim\mathcal{H}_{2\beta}^{N-1}\int_{\hat{t}\in \prod_{i=2}^N[a_i,  b_i] }\prod_{i=2}^N\sqrt{ c_i(\hat{t})}d\hat{t} u^{\frac{N-1}{\beta}}\Psi(u), \quad u\to\infty.$$
	Therefore, we conclude that
	$$\pk{\sup_{t\in E_u}X(t)>u}\sim \mathcal{H}_{2\beta}^{N-1}\int_{\hat{t}\in \prod_{i=2}^N[a_i,  b_i] }\prod_{i=2}^N\sqrt{ c_i(\hat{t})}d\hat{t} u^{\frac{N-1}{\beta}}\Psi(u), \quad u\to\infty,$$
	together with  (\ref{main})-(\ref{upper2}), establishing the claim. This completes the proof. \QED

\begin{small}

\end{small}

\bigskip

\begin{quote}
\begin{small}

\textsc{Dan Cheng}\\
School of Mathematical and Statistical Sciences\\
Arizona State University\\
900 S Palm Walk\\
Tempe, AZ 85281, U.S.A.\\
E-mail: \texttt{cheng.stats@gmail.com}

\vspace{.1in}
		
\textsc{Peng Liu}\\
Department of Actuarial Science \\
University of Lausanne\\
UNIL-Dorigny\\
1015 Lausanne, Switzerland\\
E-mail: \texttt{peng.liu1@uwaterloo.ca}

	\end{small}
\end{quote}


\begin{thebibliography}{9}
	
\bibitem{Adler00}
Adler, R. J. (2000). On excursion sets, tube formulas and maxima of random fields.
{\it Ann. Appl. Probab.} {\bf 10}, 1--74.	

\bibitem{Adler07}
Adler, R. J. and Taylor, J. E. (2007). {\it Random Fields and Geometry}. Springer, New York.






\bibitem{AzaisBook09}
Aza\"is, J.-M. and Wschebor, M. (2009).
{\it Level Sets and Extrema of Random Processes and Fields.}
John Wiley \& Sons, Inc., Hoboken, NJ.

\bibitem{Bai:2018}
Bai, L., Debicki, K., Hashorva, E. and Luo, L. (2018). On generalised Piterbarg constants. {\it Methodol. Comput. Appl. Probab.} {\bf 20}, 137--164.

\bibitem{Bickel:1973}
Bickel, P. and Rosenblatt, M. (1973). On some global measures of the deviations of density function estimates. {\it Ann. Statist.} {\bf 6}, 1071--1095.

\bibitem{ChanL06}
Chan, H. P. and Lai, T. L. (2006). Maxima of asymptotically Gaussian random fields and moderate deviation approximations to boundary crossing probabilities of sums of random variables with multidimensional indices. {\it Ann. Probab.} {\bf 34}, 80--121.


\bibitem{ChengXiaoSphere}
Cheng, D. and Xiao, Y. (2016). Excursion probability of Gaussian random fields on spheres. {\it Bernoulli}. {\bf 22}, 1113--1130.

\bibitem{Cohen12}
Cohen, S. and Lifshits, M. A. (2012). Stationary Gaussian random fields on hyperbolic spaces and on Euclidean spheres. {\it ESAIM. Probability and Statistics} {\bf 16}, 165--221.

\bibitem{CL67}
Cram\'er, H. and Leadbetter, M. R. (1967). {\it Stationary and Related Stochastic Processes:
	Sample Function Properties and Their Applications}. Wiley, New York.

\bibitem{Debicki:2016}
Debicki, K., Hashorva, E. and and Ji, L. (2016). Extremes of a class of nonhomogeneous Gaussian random fields, {\em Ann. Probab.} {\bf 44}, 984--1012.

\bibitem{Debicki:2017}
Debicki, K., Hashorva, E. and and Liu, P. (2017). Uniform tail approximation of homogenous functionals of Gaussian fields, {\em Adv. Appl. Probab.} {\bf 49}, 1037--1066.

\bibitem{Fyodorov04}
Fyodorov, Y. V. (2004). Complexity of random energy landscapes, glass transition, and absolute value of the spectral determinant of random matrices. {\it Phys. Rev. Lett.}, {\bf 92}, 240601.

\bibitem{Gelfand:2010}
Gelfand, A., Diggle, P., Guttorp, P. and Fuentes. M. (2010). {\it Handbook of Spatial Statistics}. Chapman \& Hall/CRC.

\bibitem{Gneiting13}
Gneiting, T. (2013). Strictly and non-strictly positive definite functions on spheres. {\it Bernoulli},
{\bf 19},  1327--1349.

\bibitem{Huckemann:2016}
Huckemann, S., Kim, K., Munk, A., Rehfeldt, F., Sommerfeld, M., Weickert, J. and Wollnik., C. (2016). The circular SiZer, inferred persistence of shape parameters and application to early stem cell differentiation. {\it Bernoulli}. {\bf 22}, 2113--2142.

\bibitem{HP99}
H\"usler, J. and Piterbarg, V. I. (1999). Extremes of a certain class of Gaussian processes. {\it Stoch. Process. Appl.} {\bf 83}, 257--271.

\bibitem{Istas:2005}
Istas, J. (2005). Spherical and hyperbolic fractional Brownian motion. {\em
	Electron. Comm. Probab.} {\bf 10}, 254--262.

\bibitem{Istas:2006}
Istas, J. (2006). Karhunen-Lo\`eve expansion of spherical
fractional Brownian motions. {\it Statist. Probab. Lett.} {\bf 76}, 1578--1583.


\bibitem{Jeong:2017}
Jeong, J., Jun, M. and Genton, M.G. (2017). Spherical process models for global spatial statistics. {\it Statistical Science.} {\bf 32}, 501--513.

\bibitem{Liu12}
Liu, J. (2012). Tail approximations of integrals of gaussian random fields. {\it Ann. Probab.} {\bf 40}, 1069--1104.

\bibitem{MaYang:2012}
Ma, S., Yang, L. and Carroll, R. (2012). A simultaneous confidence band for sparse longitudinal regression. {\it Statistica Sinica.} {\bf 22}, 95--122.

\bibitem{MP11}
Marinucci, D. and Peccati, G. (2011). {\it Random Fields on the Sphere. Representation, Limit Theorems and Cosmological Applications}. Cambridge University Press.

\bibitem{OhLi04}
Oh, H.-S. and Li, T.-H. (2004). Estimation of global temperature fields
from scattered observations by a spherical-wavelet-based spatially adaptive
method. {\it J. Roy. Statist. Soc. Ser. B} {\bf 66}, 221--238.

\bibitem{Piterbarg:1996}
Piterbarg, V. I. (1996). {\it Asymptotic Methods in the Theory of Gaussian Processes and Fields}. Translations of Mathematical Monographs, Vol. 148, American Mathematical Society, Providence, RI.

\bibitem{Qiao:2017}
Qiao, W. and Polonik, W. (2018). Extrema of rescaled locally stationary Gaussian fields on manifolds. {\it Bernoulli.} {\bf 24}, 1834--1859.

\bibitem{Rice:1945}
Rice, S. O. (1945). Mathematical analysis of random noise. {\it Bell System Tech. J.}, {\bf 24}, 46--156.

\bibitem{Schoenberg42}
Schoenberg, I. J. (1942). Positive definite functions on spheres. {\it Duke Math. J.},
{\bf 9},  96--108.

\bibitem{Stein:1999}
Stein, M. L. (1999). {\it Interpolation of Spatial Data: Some Theory for Kriging}. Springer, New York.

\bibitem{Stein:2007}
Stein, M. L. (2007). Spatial variation of total column Ozone on a global scale.
{\it Ann. Appl. Statist.} {\bf 1}, 191--210.

\bibitem{Sun93}
Sun, J. (1993). Tail probabilities of the maxima of Gaussian random fields. {\it Ann. Probab.} {\bf 21}, 34--71.

\bibitem{TaylorW07}
Taylor, J. E. and Worsley, K. J. (2007). Detecting sparse signals in random fields, with an application to brain mapping. {\it J. Amer. Statist. Assoc.} {\bf 102}, 913--928.

\bibitem{TaylorW08}
Taylor, J. E. and Worsley, K. J. (2008). Random fields of multivariate test statistics, with applications to shape analysis. {\it Ann. Statist.} {\bf 36}, 1--27.

\bibitem{WangYang:2007}
Wang, L. and Yang, L. (2007). Spline-backfitted kernel smoothing of nonlinear additive autoregression model. {\it Ann. Statist.} {\bf 35}, 2474--2503.

\bibitem{Yaglom61}
Yaglom, A. M. (1961). Second-order homogeneous random fields.
{\it Proc. Fourth Berkeley Symp. on Math. Statist. and Probab.},
Vol. 2, pp. 593--622. University of Calif. Press.
	
\end{thebibliography}
\end{document}